\documentclass[11pt, a4paper]{amsart}
\usepackage{color}
\def\RR{{\mathbb {R}}}
\def\I{{\mathcal I}}

\def\ep{\varepsilon}

\def\di{\displaystyle}

\def\R{{\mathbb {R}}}

\def\K{{\mathcal{K}}}
\def\L{{\mathcal{L}}}
\def\H{{\mathcal{H}}^{N-1}}

\def\pint{\operatorname {--\!\!\!\!\!\int\!\!\!\!\!--}}

\def\wg{W^{1,G}(\Omega)}
\def\wg0{W_0^{1,G}(\Omega)}

\newtheorem{teo}{Theorem}[section]
\newtheorem{lema}{Lemma}[section]
\newtheorem{defi}{Definition}[section]

\newtheorem{corol}{Corollary}[section]
\theoremstyle{definition}
\newtheorem{remark}{Remark}[section]

\begin{document}

\title[A free boundary problem in Orlicz spaces] {A minimum problem with free
boundary  in Orlicz spaces}

\author[ S. Mart\'{i}nez \& N. Wolanski]
{Sandra Mart\'{\i}nez and Noemi Wolanski}

\address{Departamento  de Matem\'atica, FCEyN
\hfill\break\indent UBA (1428) Buenos Aires, Argentina.}
\email{{ smartin@dm.uba.ar \\
wolanski@dm.uba.ar} \hfill\break\indent {\em Web-page:}
 {\tt http://mate.dm.uba.ar/$\sim$wolanski}}

\thanks{Supported by ANPCyT PICT No.
03-13719, UBA X052 and X066 and Fundaci\'on Antorchas
13900-5. N. Wolanski is a member of CONICET}


\begin{abstract}
We consider the  optimization problem of minimizing
$\int_{\Omega}G(|\nabla u|)+\lambda \chi_{\{u>0\}}\, dx$ in the
class of functions  $W^{1,G}(\Omega)$ with $u-\varphi_0\in
W_0^{1,G}(\Omega)$, for a given $\varphi_0\geq 0$ and bounded.
$W^{1,G}(\Omega)$ is the class of weakly differentiable functions
with $\int_\Omega G(|\nabla u|)\,dx<\infty$. The conditions on the
function $G$ allow for a different behavior at
 $0$ and at $\infty$. We prove that every solution $u$
is locally Lipschitz continuous, that it is a  solution to a free
boundary problem and that the free boundary, $\Omega\cap\partial\{u>0\}$, is a regular surface.
Also, we introduce the notion of weak solution to the free boundary problem
solved by the minimizers and prove the Lipschitz regularity of the weak
solutions and the $C^{1,\alpha}$ regularity of their free boundaries near
``flat'' free boundary points.
\end{abstract}

\maketitle

\section{Introduction}

In this paper we study the following  minimization problem. For
$\Omega$ a smooth bounded domain in $\R^N$   and $\varphi_0$ a
nonnegative function with $\varphi_0\in L^{\infty}(\Omega)$ and
$\int_\Omega G(|\nabla \varphi_0|)\,dx<\infty$, we consider the
problem of minimizing the functional,
\begin{equation}\label{problem}
\mathcal{J}(u)=\int_{\Omega}G(|\nabla u|)+\lambda
\chi_{\{u>0\}}\, dx
\end{equation}
 in the class of functions

$$
\mathcal{K}=\Big\{v\in L^1(\Omega):\ \int_\Omega
G(|\nabla v|)\,dx<\infty,\ v=\varphi_0\mbox{ on }\partial\Omega \Big\}.
$$

\medskip

This kind of optimization problem has been widely studied for different functions
$G$. In fact,  the first paper in which this problem was studied is \cite{AC}. The
authors considered the case $G(t)=t^2$. They proved that minimizers are weak solutions
to the free boundary problem
\begin{equation}\label{fbp}
\left\{\begin{aligned}
&\Delta u=0\quad\ \ \ \ \ \quad\quad\,\mbox{in }\{u>0\}\\
&u=0,\ |\nabla u|=\lambda\quad\mbox{on }\partial\{u>0\}
\end{aligned}\right.
\end{equation}
and proved the Lipschitz regularity of the solutions and the $C^{1,\alpha}$ regularity
of the free boundaries.

This free boundary problem appears in several applications. A very important one is
that of fluid flow. In that context, the free boundary condition is known as Bernoulli's condition.

The results of \cite{AC} have been generalized to several cases.
For instance, in \cite{ACF} the authors consider problem
\eqref{problem} for a convex function $G$ such that $ct<G'(t)<Ct$
for some  positive constants $c$ and $C$. Recently, in the article
\cite{DP} the authors considered the case $G(t)=t^p$ with
$1<p<\infty$. In these two papers only minimizers are studied.
Minimizers satisfy very good properties like nondegeneracy at the
free boundary and uniform positive density of the set $\{u=0\}$ at
free boundary points. On the other hand, the free boundary problem
\eqref{fbp} and its counterpart for different choices of functions
$G$ appears in different contexts. For instance, as limits of
singular perturbation problems of interest in combustion theory
(see for instance, \cite{BCN}, \cite{LW}). The study of weak
solutions to \eqref{fbp} also appears when considering some
optimization problems with a volume constrain (see for instance,
\cite{AAC}, \cite{ACS}, \cite{FBRW}, \cite{FBMW}, \cite{L},
\cite{MW}). Thus, the study of the regularity of weak solutions
and their free boundaries, while including the case of minimizers,
it is of a wider interest.

Thus, one of the goals of this paper is to return to the ideas of \cite{AC} and study weak solutions.
Nevertheless, our main goal is to get these results under the natural conditions on $G$ introduced by Lieberman
(see \cite{Li}) for the study of the regularity of weak solutions to the elliptic equation (possibly degenerate or
singular)
\begin{equation}\label{L}
{\mathcal L}u=\mbox{div\,}\Big(g(|\nabla u|)\frac{\nabla u}{|\nabla u|}\Big)
\end{equation}
where $g(t)=G'(t)$.

These conditions ensure that the equation \eqref{L} is equivalent to a uniformly elliptic equation in nondivergence
form with ellipticity constants independent of the solution $u$ on sets where $\nabla u\neq0$. Moreover,  these
conditions do not
imply any kind of homogeneity on the function $G$ and moreover, they allow
for a different behavior of the function $g$ when $|\nabla u|$ is close to zero or infinity.
Namely, we assume that
 $g$ satisfies
\begin{equation}\label{cond}
0<\delta\le\frac{t g'(t)}{g(t)}\le g_0\ \ \ \ \forall t>0
\end{equation}
for certain constants $\delta $ and $g_0$.

Observe that $\delta=g_0=p-1$ when $G(t)=t^p$, and conversely, if
$\delta=g_0$ then $G$ is a power. A different example consists of
a function $G$ such that $g(t)=t^a\mbox{log\,}(b t+c)$ with $a, b,
c>0$ that satisfies \eqref{cond} with $\delta=a$ and $g_0=a+1$.
Another interesting case is that of a function $G$ with  $g\in
C^1([0,\infty))$, $g(t)=c_1t^{a_1}$ for $t\le s$,
$g(t)=c_2t^{a_2}+d$ for $t\ge s$. In this case $g$ satisfies
\eqref{cond} with $\delta=\mbox{min}(a_1,a_2)$ and
$g_0=\mbox{max}(a_1,a_2)$.

Moreover, any linear combination with positive coefficients of functions satisfying \eqref{cond}
also satisfies \eqref{cond}.
Also, if $g_1$ and $g_2$ satisfy condition \eqref{cond} with constants $\delta^i$ and $g_0^i$, $i=1,2$, the function
$g=g_1g_2$ satisfies \eqref{cond} with $\delta=\delta^1+\delta^2$ and $g_0=g_0^1+g_0^2$, and the function
$g(t)=g_1\big(g_2(t)\big)$ satisfies \eqref{cond} with $\delta=\delta^1\delta^2$ and $g_0=g_0^1g_0^2$.

 This observation shows that there is a wide range of functions $G$ under
the hypothesis of this paper.

\bigskip

The main results in this article are:
\begin{teo}\label{lip}
If $g$ satisfies \eqref{cond},  there exists a  minimizer
of $\mathcal{J}$ in $\mathcal K$ and any minimizer $u$ is nonnegative and belongs to  $ C^{0,1}_{loc}(\Omega)$.
Moreover, for any domain $D\subset\subset \Omega$ containing a
free boundary point, the Lipschitz constant of $u$ in $D$ is controlled in terms of  $N, g_0, \delta,
dist(D,\partial\Omega)$ and $\lambda$.
\end{teo}

We also prove that ${\mathcal L}\,u=0$ in the set $\{u>0\}$ and that  $\{u>0\}$ has finite perimeter
locally in $\Omega$. As usual, we define the reduced boundary  by
$\partial_{red}\{u>0\}:=\{x\in \Omega\cap\partial\{u>0\}\,/\, |\nu_u(x)|=1\}$, where $\nu_u(x)$
is the unit outer normal  in the measure theoretic sense (see \cite{F}), when it exists, and $\nu_u(x)=0$
otherwise. Then, we prove that
 $\H(\partial\{u>0\}\setminus
\partial_{red}\{u>0\})=0$.

We also prove that minimizers have an the asymptotic development  near any point in their
reduced free boundary. Namely,

\begin{teo}
Let $u$ be a  minimizer, then for every
$x_0\in\partial_{red}\{u>0\}$,
  \begin{equation}\label{asim}
 u(x)=\lambda^*\langle x-x_0,\nu(x_0)\rangle^-+o(|x-x_0|)\quad \mbox{ as }
 x\rightarrow x_0
 \end{equation}
 where $\lambda^*$ is such that $g(\lambda^*)\lambda^*-G(\lambda^*)=\lambda.$ $\big($Here $\langle\cdot,\cdot\rangle$
 denotes
 the scalar product in $\RR^N$ and $v^-=-\min(v,0)$$\big)$.
\end{teo}
So that, in a weak sense minimizers satisfy,
\begin{equation}\label{ecu}
\begin{cases}
\mathcal{L}u=0 & \mbox{ in }\quad\{u>0\},\\
u=0, \quad|\nabla u|=\lambda^* & \mbox{ on
}\quad \Omega\cap\partial\{u>0\}.
\end{cases}
\end{equation}

\medskip

These results suggest that we consider weak solutions of the problem
\eqref{ecu}. We  give two different definitions of  weak
solution (Definition \ref{weak} and Definition \ref{weak2}).
Minimizers of the functional $\mathcal J$ verify both definitions of
weak solution. The main difference between these two definitions
is that for  functions satisfying Definition \ref{weak} we
have that $\H(\partial\{u>0\}\setminus
\partial_{red}\{u>0\})=0$, whereas  for functions satisfying Definition
\ref{weak2} we may have $\partial_{red}\{u>0\}=\emptyset$.
Definition \ref{weak2} is more suitable for limits of singular
perturbation problems.

 Hypotheses $(1), (2)$ and $(3)$ of
Definition \ref{weak} are similar to the ones in the definition of
weak solution to the problem studied in \cite{AC}. In our case, we
 add hypothesis $(4)$ in order to prove that weak solutions satisfying Definition \ref{weak}
also have the asymptotic
development \eqref{asim} at ${\mathcal H}^{N-1}$ almost every point of the reduced free boundary. Condition
$(4)$ is also used in the proof of  the regularity of the free
boundary. We prove the following theorem,

\begin{teo}\label{teo.main}
Let $u$ be a  weak solution.   Then, $\H$
almost every point in  the reduced free boundary
$\partial_{red}\{u>0\}$ has a neighborhood where the free boundary
 is a $C^{1,\alpha}$ surface. Moreover, if $u$ is a weak solution according to Definition \ref{weak}, the
 remainder of the free boundary has $\H -$  measure zero.
\end{teo}

We point out that we prove that, if $u$ is a weak solution, the free boundary is a $C^{1,\alpha}$ surface in a
neighborhood
of every point where $u$ has the asymptotic development \eqref{asim} for some unit vector $\nu$. We prove that
this is the case for  every point in the reduced free boundary when $u$ is a minimizer (see Theorem \ref{blow2}).
So that, if $u$ is a minimizer the reduced free boundary  is an open $C^{1,\alpha}$ surface and the remainder of
the free boundary has $\H$--measure zero.

\subsection*{Outline of the paper and technical comments}

In Section 2 we give some properties of the function $g$ and
define some spaces that we use to prove existence of minimizers.
Then, we prove some properties of solutions and subsolution of
$\mathcal{L}v=0$. We also state some real analytic properties for
functions with finite $\int_\Omega G(|\nabla u|)\,dx$  and we
prove a Cacciopoli type inequality valid for these functions. We
also prove an inequality (Theorem \ref{desilarga}) that will be
used several times in this work.

In Section 3 we prove the existence of minimizers and that they are
subsolutions of $\mathcal{L}v=0$. We also prove a maximum principle and the positivity of the minimizers.
The existence of minimizers, while standard in its form, makes strong use of the Orlicz spaces and the second inequality
in condition \eqref{cond}.

In Section 4 we prove that  any local minimizer $u$ is H\"{o}lder
continuous (Theorem \ref{cotaalpha}),   $\mathcal{L}u=0$ in
$\{u>0\}$ (Lemma \ref{solucion}) and finally we prove the local
Lipschitz continuity (Theorem \ref{Lip}). The proof of the
H\"older continuity of the minimizers is a key step in our
analysis. Although, the proof follows closely the one for the case of the
$p-$laplacian (\cite{DP}), here we have to use all the properties
of the function $G$  which mainly come into play through the  inequality in
Theorem \ref{desilarga}.

In Section 5 we prove that minimizers satisfy a nondegeneracy
property near the free boundary $\Omega\cap
\partial\{u>0\}$. We also prove that the sets
$\{u>0\}$ and $\{u=0\}$ have locally uniform positive density at the free boundary (Theorem
\ref{densprop}). In this theorem we make strong use of
the properties of $G$ and the corresponding Orlicz space.

In Section 6 we prove that the free boundary  has Hausdorff dimension $N-1$ and  we
obtain  a representation
theorem for minimizers (Theorem \ref{repteo}). This implies that
$\{u>0\}$ has locally finite perimeter  in $\Omega$.
Finally we prove that $\H(\partial\{u>0\}\setminus\partial_{\rm
red}\{u>0\})=0.$

In Section 7 we give some properties of blow up sequences of
minimizers. We prove that any limit of a blow up sequence of minimizers is again a
minimizer (Lemma \ref{minblow}) and we finally prove the asymptotic
development of minimizers at every point in their reduced free boundary  (Theorem
\ref{blow2}).

In Section 8 we give the definition of weak solution (Definition
\ref{weak} and Definition \ref{weak2}). We show that most of the
properties that we proved for minimizers also hold for weak
solutions according to Definition \ref{weak},  and we mention the
differences between the two definitions (Remark \ref{repweak} and
Remark \ref{diferencia}).

In Section 9 we prove the regularity  of the free boundary
of weak solutions near ``flat'' free boundary points (Theorem \ref{regfin})
and then, we deduce the regularity of the free boundary of weak solutions near
almost every point in their reduced free boundary and, in the case of minimizers, the  regularity
of the whole reduced free boundary  (Theorem \ref{teo.regularity}). While most of the steps of the proof of
the regularity of the free boundary of weak solutions are very similar to the corresponding ones for minimizers
 in the uniformly elliptic case considered in \cite{ACF} and
in the case $G(t)=t^p$ considered in \cite{DP}, there are some steps that need a new proof since weak solutions do not
verify the locally uniform positive density of $\{u=0\}$ at the free boundary (See Lemmas \ref{9.1} and \ref{flat3} and
Theorem \ref{regfin}).

\bigskip


\section{Properties of the function $G$}

In this section we state and prove some properties of the function
$G$ and its derivative $g$ that are used throughout the paper. We
also state some real analytic properties for functions with finite
$\int_\Omega G(|\nabla u|)\,dx$  like a form of Poincar\'e
Inequality, a Cacciopoli type inequality, the H\"older continuity
of functions in a kind of Morrey type space, properties of weak
solutions to $\L u=0$ and a comparison principle for sub and
supersolutions. We also prove an important inequality (Theorem
\ref{desilarga}). All these properties will be thoroughly used
throughout the paper. Some of them have been proved in \cite{Li}.
We only write down the proof of statements not contained in
\cite{Li}.

\begin{lema}\label{prop}
The function $g$ satisfies the following properties,
\begin{enumerate}
\item[(g1)] $\di \min\{s^{\delta},s^{g_0}\} g(t)\leq g(st)\leq
\max\{s^{\delta},s^{g_0}\} g(t)$
\smallskip
 \item[(g2)] $G$ is convex and $C^2$\item[(g3)]
$\di\frac{t g(t)}{1+g_0}\leq G(t)\leq t g(t) \quad \forall\ t\geq
0$

\end{enumerate}

\end{lema}
\begin{proof}
For the proofs of (g1)--(g3) see \cite{Li}.
\end{proof}
\begin{remark} By (g1) and (g3) we have a similar inequality for
$G$,
\begin{enumerate}
\item[(G1)]  $\di\min\{s^{\delta+1},s^{g_0+1}\}
\frac{G(t)}{1+g_0}\leq G(st)\leq
(1+g_0)\max\{s^{\delta+1},s^{g_0+1}\} G(t)$
\end{enumerate}
and, then using the convexity of $G$ and this last inequality we
have, \begin{enumerate} \item[(G2)] $\displaystyle G(a+b)\leq
{2^{g_0}}{(1+g_0)} (G(a)+G(b))\ \forall\ a,b>0$.
\end{enumerate}
\end{remark}

As $g$ is strictly increasing we can define $g^{-1}$. Now we prove
that $g^{-1}$ satisfies a condition similar to \eqref{cond}. That
is,

\begin{lema}
The function $g^{-1}$ satisfies the inequalities
\begin{equation}\label{condinv}\frac{1}{g_0} \leq \frac{t {(g^{-1})}'(t)}{g^{-1}(t)}\leq
\frac{1}{\delta} \ \ \ \ \forall t>0.\end{equation}
Moreover, $g^{-1}$ satisfies,
\begin{equation}\label{gmenos1}\tag{$\widetilde{g}   1$}
 \di \min\{s^{1/\delta},s^{1/g_0}\}
g^{-1}(t)\leq g^{-1}(st)\leq \max\{s^{1/\delta},s^{1/g_0}\}
g^{-1}(t)\end{equation} and if $\widetilde{G}$ is such that
${\widetilde{G}}'(t)=g^{-1}(t)$ then,
\begin{equation}\label{gmenos2}\tag{$\widetilde{g} 2$} \di\frac{\delta t
g^{-1}(t)}{1+\delta}\leq \widetilde{G}(t)\leq t g^{-1}(t) \quad
\forall\ t\geq 0\end{equation}
\begin{equation}\label{Gmono1}\tag{$\widetilde{G}1$}
\di\frac{(1+\delta)}{\delta}\min\{s^{1+1/\delta},s^{1+1/g_0}\}\widetilde{G}(t)
\leq \widetilde{G}(st)\leq
\frac{\delta}{1+\delta}\max\{s^{1+1/\delta},s^{1+1/g_0}\}
\widetilde{G}(t)\end{equation}
\begin{equation}\label{gmenos3}\tag{$\widetilde{g} 3$} \di a b\leq \ep G(a)+
C(\ep)\widetilde{G}(b)\quad \forall \ a,b>0\mbox{ and }\ep>0\mbox{ small}\end{equation}
\begin{equation}\label{gmenos4}\tag{$\widetilde{g} 4$} \di\widetilde{G}(g(t))\leq
g_0 G(t)\end{equation}
\end{lema}
\begin{proof}
Let $s=g^{-1}(t)$, then
$$\frac{t {(g^{-1})}'(t)}{g^{-1}(t)}=\frac{
g(s)}{g'(s)s}$$ and using \eqref{cond} we have the desired
inequalities.

Now \eqref{gmenos1} follows by property (g1) applied to $g^{-1}$,
and ($\widetilde{g}   2$) by property (g3). ($\widetilde{G}1$)
follows by $\widetilde{g}_1$ and $\widetilde{g}_2$.
\\
By Young's inequality  we have that $ab\leq G(a)+\widetilde{G}(b)$
and then, for $0<\ep'<1$ such that $\ep=(1+g_0)\ep'^{(1+\delta)}$,
$$\ep' a \frac{b}{\ep'}\leq G(\ep'
a)+\widetilde{G}\Big(\frac{b}{\ep'}\Big)\leq \ep\, G(a)+ C(\ep)
\widetilde{G}(b).$$ In the las inequality we have used (${G}
1$) and ($\widetilde{G}   1$). Thus ($\widetilde{g}   3$) follows.\\
As $g$ is strictly increasing we have that
$\widetilde{G}(g(t))+G(t)=tg(t)$ (see
 equation (5), Section 8.2 in \cite{Ad}) and applying
 (g3),
  we get
 $$\widetilde{G}(g(t))=t g(t)-G(t)\leq g_0 G(t).$$ Thus, ($\widetilde{g}   4$)
 follows.
\end{proof}

\medskip

In order to prove the existence of minimizers we will use some compact
embedding results. To this end, we have to define some Orlicz and
Orlicz-Sobolev spaces. We recall that the functional
$$\|u\|_G=\inf\Big\{k>0:\int_{\Omega} G\Big(\frac{|u(x)|}{k}\Big)\,
dx \leq 1\Big\}$$ is a norm in the Orlicz space $L^{G}(\Omega)$
which is the linear hull of the Orlicz class
$$K_G(\Omega)=\Big\{u \mbox{ measurable }:\ \int_{\Omega}
G(|u|)\, dx<\infty\Big\},$$ observe that this set is convex,
since $G$ is also convex (property (g2)). The Orlicz-Sobolev space
$W^{1,G}(\Omega)$ consists of those functions in $L^{G}(\Omega)$
whose distributional derivatives $\nabla u$ also belong to
$L^G(\Omega)$. And we have that
$\|u\|_{W^{1,G}}=\max\{\|u\|_G,\|\nabla u\|_G\}$ is a norm for
this space.

\begin{lema}\label{equi}
There exists a constant $C=C(g_0,\delta)$ such that,
$$\|u\|_G \leq C \max\Big\{\Big(\int_{\Omega} G(|u|)\, dx\Big)^{1/(\delta+1)},
\Big(\int_{\Omega} G(|u|)\, dx\Big)^{1/{(g_0+1)}}\Big\}$$
\end{lema}
\begin{proof}

\end{proof}
If $\di\int_{\Omega} G(|u|)\, dx= 0$ then $u=0$ a.e and the result
follows. If $\int_{\Omega} G(|u|)\, dx\neq 0$, take
$k=\max\Big\{\Big(2(1+g_0)\int_{\Omega} G(|u|)\, dx\Big)
^{1/{(\delta+1)}}, \Big(2(1+g_0)\int_{\Omega} G(|u|)\, dx\Big)
^{1/{(g_0+1)}}\Big\}$,  by (G1) we have,
$$\int_{\Omega} G\Big(\frac{|u|}{k}\Big)\, dx\leq (1+g_0) \max\Big\{\frac{1}{k^{\delta+1}},
\frac{1}{k^{g_0+1}} \Big\}\int_{\Omega} G(|u|)\, dx\leq 1$$
therefore $\|u\|_G \leq k$ and the result follows.

\begin{teo}\label{ref}
$L^{\widetilde{G}}(\Omega)$ is the dual of $L^G(\Omega)$.
Moreover, $L^G(\Omega)$ and $W^{1,G}(\Omega)$ are reflexive.
\end{teo}
\begin{proof}
As $G$  satisfies property  $(G1)$ and $\widetilde{G}$ property
$(\widetilde{G} 1)$, we have that both pairs $(G,\Omega)$ and
$(\widetilde{G},\Omega)$ are $\Delta-$ regular (see 8.7 in
\cite{Ad}). Therefore we are in the hypothesis of Theorem 8.19 and
Theorem 8.28 at \cite{Ad}, and the result follows.
\end{proof}
\begin{teo}\label{imb}
 $L^G(\Omega)\hookrightarrow L^{1+\delta}(\Omega)$
continuously.
\end{teo}
\begin{proof}
By theorem 8.12 of \cite{Ad} we only have to prove that $G$
dominates $t^{1+\delta}$ near infinity. That is, there exits
constants $k$, $t_0$ such that $t^{1+\delta}\leq G(kt)\quad
\forall t\ \geq t_0$. But this is true by property (G1). So the
result follows.
\end{proof}

\medskip

The following result is a Poincar\'e type inequality.

\begin{lema}\label{poinc} If $u\in W^{1,1}(\Omega)$ with $u=0$ on $\partial\Omega$ and
 $\int_{\Omega} G(|\nabla u|)\,
dx$ is finite, then
$$\int_{\Omega} G\Big(\frac{|u|}{R}\Big)\, dx \leq
\int_{\Omega} G(|\nabla u|)\, dx \quad \mbox{ for
}R=\mbox{diam\,}\Omega.$$
\end{lema}
\begin{proof}
See Lemma 2.2 of \cite{Li}.
\end{proof}

\medskip

Now we state a generalization of Morrey's Theorem. Let
$$
[u]_{0,\alpha,\Omega}= \sup_{\stackrel{x,y\in \Omega}{x\neq y}}
\frac{|u(x)-u(y)|}{|x-y|^{\alpha}}.
$$

We have the following
result,

\begin{lema}\label{calfa}
Let $u\in L^{\infty}(\Omega)$  such that for some $0<\alpha<1$ and $r_0>0$,
$$\int_{B_r} G(|\nabla
u|)\, dx \leq C r^{N+\alpha-1}\quad \forall\ 0<r\leq r_0
$$
with $B_r\subset \Omega$. Then, $u\in C^{\alpha}(\Omega)$ and
there exists a constant $C_1=C_1(C,\alpha,N,g_0,G(1))$ such that
$[u]_{0,\alpha,\Omega}\leq  C_1$.

\end{lema}
\begin{proof}
The proof of this lemma is included in the proof of Theorem 1.7
(pag. 346) in \cite{Li}.
\end{proof}

Now, we will give some properties of subsolutions and solutions of
 $\mathcal{L}v=\mbox{div}(A(\nabla v))=0$, where $A(p)=g(|p|)
\frac{ p}{| p|}$. First, let us observe that if
$a_{ij}=\frac{\partial A_i}{\partial p_j}$ by using
\eqref{cond}, we get
\begin{equation}\label{aij}
\min\{\delta,1\} \frac{ g(|p|)}{| p|}|\xi|^2\leq a_{ij} \xi_i
\xi_j \leq \max\{g_0,1\} \frac{ g(|p|)}{| p|}|\xi|^2,
\end{equation}
which means that the equation $\mathcal{L}v=0$ is uniformly
elliptic for $\frac{ g(|p|)}{| p|}$ bounded and bounded away from
zero.

\medskip

The next lemma is a Cacciopoli type inequality for subsolutions of $\mathcal{L} v=0.$

\begin{lema}\label{caccio}
Let $v$ be a nonnegative weak subsolution of $\mathcal{L} v=0.$
That is,
\begin{equation}\label{subs}0\geq \int_{\Omega} g(|\nabla v|) \frac{\nabla
v}{|\nabla v|} \nabla \phi\, dx\quad \forall\ \phi \in
C_0^{\infty}(\Omega) \mbox{ such that } \phi\geq 0.\end{equation}
Then, there exists $C=C(N,\delta, g_0)>0$ such that
$$\int_{B_r} G(|\nabla v|)\, dx \leq C \int_{B_{\frac{3}{2}r}}
G\Big(\frac{|v|}{r}\Big)\, dx $$ for all $r>0$, such that
$B_{\frac{3}{2}r}\subset\Omega$.
\end{lema}
\begin{proof}
Let $\phi=v \eta^{g_0+1}$, where $0\le \eta\in
C^1_0(B_{\frac{3}{2}r})$, with $|\nabla \eta|\leq \frac{C}{r}$,
$\eta\leq 1$, $\eta\equiv 1 \mbox{ in } B_r$. Then, $\nabla
\phi=\eta^{g_0+1}\nabla v + v \nabla \eta (g_0+1) \eta^{g_0}$ and
replacing in \eqref{subs} we have,
$$0\geq \int_{B_{\frac{3}{2}r}} g(|\nabla v|) |\nabla v| \eta^{g_0+1} \, dx +(g_0+1)
\int_{B_{\frac{3}{2}r}} g(|\nabla v|) \frac{\nabla v}{|\nabla v|}
\nabla \eta\ v\,  \eta^{g_0} \, dx.
$$
Then,
$$ \int_{B_{\frac{3}{2}r}} g(|\nabla v|) |\nabla v| \eta^{g_0+1} \, dx\leq
(g_0+1) \int_{B_{\frac{3}{2}r}} g(|\nabla v|) |\nabla \eta| |v|
\eta^{g_0} \, dx,
$$
By property $(\widetilde{g}3)$ we have,
$$g(|\nabla v|) |\nabla \eta| |v|  \eta^{g_0}\leq \ep
\widetilde{G}(g(|\nabla v|) \eta^{g_0})+C(\ep) G(|\nabla \eta|
|v|).$$  Then, by property $(\widetilde{G}1)$ and as $\eta \leq
1$, we have,
$$\widetilde{G}(g(|\nabla v|) \eta^{g_0})\leq C
\eta^{g_0\big(1+\frac{1}{g_0}\big)} \widetilde{G}(g(|\nabla
v|))\leq C \eta^{1+g_0} G(|\nabla v|),$$ where the last inequality
holds by $(\widetilde{g}4)$. Summing up, and using property (g3),
we obtain
$$\int_{B_{\frac{3}{2}r}}G(|\nabla v|)  \eta^{g_0+1}\, dx\leq
C\ep \int_{B_{\frac{3}{2}r}} G(|\nabla v|) \eta^{g_0+1}\,
dx+C(\ep)\int_{B_{\frac{3}{2}r}} G(|\nabla \eta| |v|)\, dx,$$ and
if we take $\ep$ small and use the bound for $|\nabla\eta|$ we
have,
$$\int_{B_{\frac{3}{2}r}}G(|\nabla v|)  \eta^{g_0+1}\,
dx\leq C \int_{B_{\frac{3}{2}r}} G(|\nabla \eta| |v|)\, dx \leq C
\int_{B_{\frac{3}{2}r}} G\Big(\frac{|v|}{r}\Big)\, dx.$$
Finally, if we use that $\eta\equiv 1$ in $B_r$ the result
follows.

\end{proof}

\begin{lema}\label{propsol}
Let $v$ be a weak solution of $\mathcal{L}v=0$, that is
$$\int_{\Omega} g(|\nabla u|) \frac{\nabla v}{|\nabla v|}
\nabla \phi\, dx=0\quad \forall\ \phi \in C_0^{\infty}(\Omega).$$
Then $v\in C^{1,\alpha}(\Omega)$. Moreover, there exists
$C=C(N,\delta,g_0)>0$  such that for every ball $B_r\subset\Omega$,
\begin{equation}\tag{1}
\sup_{B_{r/2}} G(|\nabla v|)\leq \frac{C}{r^N}
\int_{B_{\frac{2}{3}r}} G(|\nabla v|)\, dx\end{equation}
\begin{equation} \tag{$2$}
\sup_{B_{r/2}} |\nabla v|\leq \frac{C}{r} \sup_{B_{r}} | v|
\end{equation}
$\mbox{ For every } \beta\in(0,N), \mbox{ there exists  } C=C(N,\beta,\delta,g_0,\|v\|_{L^\infty(\frac23r)})>0
\mbox{ such that,}$
\begin{equation} \tag{$3$}
 \int_{B_{r/2}} G(|\nabla v|)\leq C r^{\beta}.\end{equation}
\end{lema}
\begin{proof}
For the proof of (1) see Lemma 5.1 of \cite{Li} and for the proof
of (3) see (5.9) page 346  of \cite{Li}. Let us prove (2). By
using (1) and then Lemma \ref{caccio} we have,
$$\sup_{B_{r/2}} G(|\nabla v|)\leq \frac{C}{r^N}
\int_{B_{\frac{2}{3}r}} G(|\nabla v|)\, dx \leq \frac{C}{r^N}
\int_{B_{r}} G\Big( \frac{|v|}{r}\Big)\, dx\leq
G\Big(\frac{C}{r}\|v\|_{L^\infty(B_{r})}\Big).$$ Then
$$|\nabla v (y_0)|\leq \frac{C}{r}
\|v\|_{L^\infty(B_{r})}\quad \forall y_0\in B_{r/2}$$ and the
result follows.
\end{proof}

\begin{lema}\label{comp}
 Let $U$ be an open subset,  $u$  a weak subsolution and $w$ a weak supersolution of  $\mathcal{L}u=0$ in
 $U$.
If $w\geq u$ on $\partial U$ then, $w\geq u$ in $U$. If $w$ is a solution to $\mathcal{L}w=0$ and $w=u$ on
$\partial U$ then, $w$ is uniquely determined.
\end{lema}
\begin{proof}
\begin{align*}
0 \geq & \int_U \Big( g(|\nabla u|)\frac{\nabla u}{|\nabla
u|}-g(|\nabla w|)\frac{\nabla w}{|\nabla w|}\Big). \nabla
(u-w)^+\, dx
\\=
& \int_{U\cap \{u>w\}} \Big( g(|\nabla u|)\frac{\nabla u}{|\nabla
u|}-g(|\nabla w|)\frac{\nabla w}{|\nabla w|}\Big). \nabla
(u-w)\, dx
\\
=& \int_{U\cap \{u>w\}} \int_0^1 a_{ij}(\nabla u+(1-t) (\nabla
w-\nabla u)) (u_{x_i}-w_{x_i}) (u_{x_j}-w_{x_j})\, dt\, dx
\end{align*}
And using \eqref{aij} we have that the right hand side is grater
an equal than
\begin{align*}
C  \int_{U\cap \{u>w\}} \int_0^1  F(|\nabla u+(1-t) (\nabla
w-\nabla u)|) |\nabla w-\nabla u|^2 \, dt \, dx ,
\end{align*} where $F(t)=\frac{g(t)}{t}$.
Now, we take the following subsets of $U$
$$S_1=\{x\in U: |\nabla u-\nabla w|\leq 2 |\nabla u|\},\quad
S_2=\{x\in U: |\nabla u-\nabla w| > 2 |\nabla u|\}$$ Then $S_1\cup
S_2= U$ and
\begin{align}\label{desig}&\frac{1}{2}|\nabla u|\leq |\nabla u+(1-t)
(\nabla w-\nabla u)|\leq 3|\nabla u| &\quad \mbox{ in } S_1 \mbox{
for } t\geq \frac{3}{4}\\ \label{desig2}& \frac{1}{4}|\nabla
u-\nabla w|\leq |\nabla u+(1-t) (\nabla w-\nabla u)|\leq 3|\nabla
u-\nabla w| & \quad \mbox{ in } S_2 \mbox{ for } t\leq
\frac{1}{4}.
\end{align}
In $S_1$, and for $t\geq 3/4$ we have using \eqref{desig}, that
$$F(|\nabla u+(1-t)(\nabla w-\nabla u)|)=\frac{g(|\nabla u+(1-t)
(\nabla w-\nabla u)|)}{|\nabla u+(1-t) (\nabla w-\nabla u)|}\geq
\frac{g(\frac{1}{2}|\nabla u|)}{3|\nabla u|}\geq
\frac{1}{2^{g_0}3}F(|\nabla u|)$$ where in the last inequality
we have used (g1).

In $S_2$, and for $t\leq 1/4$ we have using (g3) and then
\eqref{desig2} that,
\begin{align*}F(|\nabla u+(1-t)(\nabla w-\nabla
u)|)|\nabla u-\nabla w|^2&\geq\frac{G(|\nabla u+(1-t) (\nabla
w-\nabla u)|)}{|\nabla u+(1-t) (\nabla w-\nabla u)|^2}|\nabla
u-\nabla w|^2\\& \geq \frac{G(\frac{1}{4}|\nabla u-\nabla
w|)}{9|\nabla u-\nabla w|^2}|\nabla u-\nabla w|^2\\&\geq
\frac{G(|\nabla u-\nabla w|)}{4^{g_0+1}9(1+g_0)}\end{align*}
where in the last inequality we have used (G1).

Therefore, we have that
$$0\geq C  \Big(\int_{S_1}   F(|\nabla u|) |\nabla (u-w)^+|^2 \, dx  +
\int_{S_2}  G(|\nabla (u-w)^+|) \, dx \Big).$$ Hence $\nabla
(u-w)^+=0$ in $S_2$
 and $\nabla (u-w)^+=0$,  or $F(|\nabla
u|)=0$ in $S_1$ in which case  $\nabla u=0$ and, by the definition of $S_1$,
this implies that $\nabla (u-w)=0$ in $S_1$. Therefore, $\nabla (u-w)^+=0$
in $U$, then $(u-w)^+=0$, which implies $u\leq w$.
\end{proof}

The following inequality will be a key tool in the proof of  the
H\"{o}lder continuity of minimizers. As an observation, we mention
that the following result is a generalization of well known
integral inequalities for the $p$- Laplacian (see, for example,
pag.4 in \cite{DP}). Here the difference is that we obtain a
unique inequality for any $\delta$ and $g_0$, (for the
$p$-Laplacian the inequalities were separated in two cases  $p\geq
2$ and
 $1<p<2$).

\begin{teo}\label{desilarga}
Let $u\in W^{1,G}(\Omega)$, $B_r\subset\subset \Omega$ and  $v$ be
a solution of
$$
\mathcal{L} v=0\quad\mbox{in }B_r, \qquad v-u\in W_0^{1,G}(B_r).
$$
then
$$
\int_{B_r} (G(|\nabla u|)-G(|\nabla v|))\, dx \geq  C\Big(\int_{
A_2} G(|\nabla u-\nabla v|)  \, dx+\int_{A_1}F(|\nabla u|) |\nabla
u-\nabla v|^2 \, dx\Big),
$$
where $F(t)=g(t)/t,$
$$A_1=\{x\in {B_r}: |\nabla u-\nabla v|\leq 2 |\nabla
u|\}\quad\mbox{ and }\quad A_2=\{x\in {B_r}: |\nabla u-\nabla v| >
2 |\nabla u|\}$$ and $C=C(g_0,\delta)$.
\end{teo}
\begin{proof}

Let $u^s=s u + (1-s) v$. Using the integral form of the mean value
theorem and the fact that $v$ is an $\mathcal{L}$-- solution, we
have,
\begin{align*}
 \I:&=\int_{B_r} (G(|\nabla u|)-G(|\nabla
v|))\, dx = \int_0^1  \int_{{B_r}} g(|\nabla u^s|)\frac{\nabla
u^s}{|\nabla u^s|} . \nabla (u-v)\, dx\, ds
\\&= \int_0^1 \frac{1}{s} \int_{B_r} \Big(
g(|\nabla u^s|)\frac{\nabla u^s}{|\nabla u^s|}-g(|\nabla
v|)\frac{\nabla v}{|\nabla v|}\Big). \nabla (u^s-v)\, dx\, ds
\\&= \int_0^1 \frac{1}{s} \int_{B_r} \int_0^1 a_{ij}(\nabla u^s+(1-t)
(\nabla v-\nabla u^s)) (u^s_{x_i}-v_{x_i}) (u^s_{x_j}-v_{x_j})\,
dt\, dx\, ds.
\end{align*}
And, by \eqref{aij} we have that the right hand side is grater
than or equal to
\begin{align*}
C \int_0^1 \frac{1}{s} \int_{B_r} \int_0^1  F(|\nabla u^s+(1-t)
(\nabla v-\nabla u^s)|) |\nabla v-\nabla u^s|^2 \, dt \, dx \, ds.
\end{align*} where $F$ was defined in Lemma \ref{comp} and
$C=C(\delta)$.

Now, we take the following subsets of ${B_r}$
$$S_1=\{x\in {B_r}: |\nabla u^s-\nabla v|\leq 2 |\nabla u^s|\},\quad
S_2=\{x\in {B_r}: |\nabla u^s-\nabla v| > 2 |\nabla u^s|\}$$ Then
$S_1\cup S_2= {B_r}$ and
\begin{align}\label{desi}&\frac{1}{2}|\nabla u^s|\leq |\nabla u^s+(1-t)
(\nabla v-\nabla u^s)|\leq 3|\nabla u^s| & \mbox{ on } S_1 \mbox{
for } t\geq \frac{3}{4}\\ \label{desi2}& \frac{1}{4}|\nabla
u^s-\nabla v|\leq |\nabla u^s+(1-t) (\nabla v-\nabla u^s)|\leq
3|\nabla u^s-\nabla v| &  \mbox{ on } S_2 \mbox{ for } t\leq
\frac{1}{4}.
\end{align}
Proceeding as in Lemma \ref{comp}, we get
$$\begin{aligned}F(|\nabla u^s+(1-t)(\nabla v-\nabla
u^s)|)
\geq \frac{1}{2^{g_0}3}F(|\nabla u^s|)\end{aligned}$$ in $S_1$
and
$$\begin{aligned}
F(|\nabla u^s+(1-t)&(\nabla v-\nabla u^s)|)|\nabla u^s-\nabla v|^2
\geq \frac{G(|\nabla u^s-\nabla
v|)}{4^{g_0+1}9(1+g_0)}\end{aligned}
$$
in $S_2$. 

Therefore, we have that
$$\I\geq C  \Big(\int_0^1\frac{1}{s}\int_{S_1}   F(|\nabla u^s|) |\nabla v-\nabla u^s|^2 \, dx \, ds +
\int_0^1 \frac{1}{s}\int_{S_2} G(|\nabla u^s-\nabla v|) \, dx \,
ds\Big)$$ Now, let
$$A_1=\{x\in {B_r}: |\nabla u-\nabla v|\leq 2 |\nabla u|\},\quad
A_2=\{x\in {B_r}: |\nabla u-\nabla v| > 2 |\nabla u|\},$$ then
${B_r}=A_1\cup A_2$, and
\begin{align}\label{desi3}&\frac{1}{2}|\nabla u|\leq |\nabla u^s|\leq 3|\nabla u| &\quad \mbox{ on } A_1
\mbox{ for } s\geq \frac{3}{4}\\ \label{desi4}& \frac{1}{4}|\nabla
u-\nabla v|\leq |\nabla u^s|\leq 3|\nabla u-\nabla v| & \quad
\mbox{ on } A_2 \mbox{ for } s\leq \frac{1}{4}.
\end{align}
Therefore \begin{align*} \I&\geq  C\Big(
\int_0^{1/4}\frac{1}{s}\int_{S_1\cap A_2} F(|\nabla u^s|) |\nabla
v-\nabla u^s|^2 \, dx \, ds \\&\quad+
\int_{3/4}^{1}\frac{1}{s}\int_{S_1\cap A_1} F(|\nabla u^s|)
|\nabla v-\nabla u^s|^2 \, dx \, ds\Big.
\Big.\\&\quad+\int_0^{1/4} \frac{1}{s}\int_{S_2\cap A_2} G(|\nabla
u^s-\nabla v|) \, dx \, ds &\\&\quad+\int_{3/4}^{1}
\frac{1}{s}\int_{S_2\cap A_1} G(|\nabla u^s-\nabla v|) \, dx \,
ds\Big)=I+II+III+IV.
\end{align*}
Let us estimate these four terms,

In $S_1\cap A_2$, for $s\leq 1/4$ we have by \eqref{desi4} and
(g1),
 that
$$F(|\nabla u^s|)\geq \frac{1}{4^{g_0}3}F(|\nabla u-\nabla
v|).$$ Therefore,
\begin{align*}I&\geq C \int_0^{1/4}\frac{1}{s}\int_{S_1\cap A_2}
F(|\nabla u-\nabla v|) |\nabla v-\nabla u^s|^2 \, dx \, ds \\&=C
\int_0^{1/4}s\int_{S_1\cap A_2} F(|\nabla u-\nabla v|) |\nabla
v-\nabla u|^2 \, dx \, ds\\& \geq C\int_0^{1/4}s\int_{S_1\cap A_2}
G(|\nabla u-\nabla v|)  \, dx \, ds\end{align*} where in the last
inequality we are using (g3).

In $S_1\cap A_1$, for $s\geq 3/4$ we have by \eqref{desi3} and
(g1), that
$$F(|\nabla u^s|)\geq \frac{1}{2^{g_0}3}F(|\nabla u|).$$
Therefore,
\begin{align*}II&\geq C \int_{3/4}^{1} s\int_{S_1\cap
A_1} F(|\nabla u|) |\nabla v-\nabla u|^2 \, dx \, ds\\&\geq C
\int_{3/4}^{1} \int_{S_1\cap A_1} F(|\nabla u|) |\nabla v-\nabla
u|^2 \, dx \, ds.\end{align*}

In $S_2\cap A_2$, for $s\leq 1/4$ we have by definition of $S_2$,
 by \eqref{desi4} and (G1), that
$$G(|\nabla u^s-\nabla v|)\geq \frac{1}{2^{g_0+1}(g_0+1)} G(|\nabla u-\nabla v|),$$ therefore
\begin{align*}III&\geq C \int_0^{1/4}\frac{1}{s}\int_{S_2\cap A_2} G(|\nabla
u-\nabla v|) \, dx \, ds\\ &\geq C \int_0^{1/4}s\int_{S_2\cap A_2}
G(|\nabla u-\nabla v|) \, dx \, ds.\end{align*}

In $S_2\cap A_1$, for $s\geq 3/4$ we have, by definition of $S_2$
and by \eqref{desi3} \begin{equation}\label{desi5}|\nabla
u^s-\nabla v|>2|\nabla u^s|\geq |\nabla u|\end{equation}

By (g3),  using \eqref{desi5} and the definition of $A_1$ we have,
\begin{align*}G(|\nabla u^s-\nabla v|)&\geq \frac1{g_0+1} g(|\nabla u^s-\nabla
v|)|\nabla u^s-\nabla v|\geq \frac1{g_0+1} g(|\nabla u|)|\nabla
u^s-\nabla v|\\& = \frac1{g_0+1} F(|\nabla u|) s |\nabla u-\nabla
v| |\nabla u|\geq  \frac s {2(g_0+1)}F(|\nabla u|) |\nabla
u-\nabla v|^2.\end{align*} Therefore,
$$ IV \geq C \int_{3/4}^{1}\int_{S_2\cap
A_1}F(|\nabla u|) |\nabla u-\nabla v|^2 \, dx\, ds.$$

If we sum $I+III$, we obtain \begin{align*}I+III&\geq
\int_0^{1/4}Cs\Big(\int_{S_1\cap A_2} G(|\nabla u-\nabla v|)  \,
dx + \int_{S_2\cap A_2} G(|\nabla u-\nabla v|) \, dx \,
ds\Big)\\&=C\int_0^{1/4}s\int_{ A_2} G(|\nabla u-\nabla v|)  \, dx
\, ds=C\int_{ A_2} G(|\nabla u-\nabla v|)  \, dx
\end{align*}
and if we sum  $II+IV$, we obtain
\begin{align*}II+IV&\geq
C\int_1^{3/4} \Big(\int_{S_1\cap A_1}F(|\nabla u|) |\nabla
u-\nabla v|^2 \, dx \Big.\\&\quad+\int_{S_2\cap A_1}F(|\nabla u|)
|\nabla u-\nabla v|^2 \, dx\Big)\, ds=C\int_{A_1}F(|\nabla u|)
|\nabla u-\nabla v|^2 \, dx.
\end{align*}
Therefore,
\begin{equation}\label{a1a2}
\I\geq  C\Big(\int_{ A_2} G(|\nabla u-\nabla v|)  \,
dx+\int_{A_1}F(|\nabla u|) |\nabla u-\nabla v|^2 \, dx\Big),
\end{equation}
where $C=C( g_0,\delta)$.

\end{proof}

In Section 4 we will need an explicit family of subsolutions and
supersolutions in an annulus. We state here the required lemma.

\begin{lema}\label{exp}
Let $w_{\mu}=\ep e^{-\mu|x|^2}$, for $\ep>0$, $r_1>r_2>0$ then
there exists $\mu>0$ such that
\begin{align*}&\mathcal{L}w_{\mu}>0
 \mbox{ in }B_{r_1}-B_{r_2}\end{align*}
 and $\mu$ depends only on $r_2, g_0,\delta$ and $N$.
\end{lema}
\begin{proof}
First, note that $$\mathcal{L}w=\frac{g(|\nabla w|)}{|\nabla w|^3}
\Big\{\Big( \frac{g'(|\nabla w|)}{g(|\nabla w|)}|\nabla w|-1
\Big)\sum_{i,j} w_{x_i}w_{x_j}w_{x_ix_j}+\triangle w|\nabla w|^2
\Big\}.$$ Computing, we have
\begin{align}\label{gradw}
w_{x_i}=-2\ep\mu x_i e^{-\mu|x|^2}& &w_{x_i x_j}=\ep(4\mu^2 x_i
x_j-2\mu \delta_{ij}) e^{-\mu|x|^2}& &|\nabla w|=2\ep\mu
|x|e^{-\mu|x|^2},
\end{align}
therefore using \eqref{gradw} and \eqref{cond} we obtain,
\begin{align*}\di e^{3\mu|x|^2}\mathcal{L}w&=\ep^3\di\frac{g(|\nabla w|)}{|\nabla w|^3} \Big\{\Big(
\frac{g'(|\nabla w|)}{g(|\nabla w|)}|\nabla w|-1 \Big)(16
\mu^4|x|^4-8\mu^3 |x|^2)+ (4\mu^2|x|^2-2\mu
N)4\mu^2|x|^2\Big\}\\&\di=\ep^3 \frac{g(|\nabla w|)}{|\nabla w|^3}
4\mu^3|x|^2\Big\{\Big( \frac{g'(|\nabla w|)}{g(|\nabla w|)}|\nabla
w|-1 \Big)(4 \mu|x|^2-2)+ (4\mu|x|^2-2 N)\Big\}\\\di&=\ep^3
\frac{g(|\nabla w|)}{|\nabla w|^3} 4\mu^3|x|^2\Big\{\Big(
\frac{g'(|\nabla w|)}{g(|\nabla w|)}|\nabla w| \Big)4
\mu|x|^2-\Big(\frac{g'(|\nabla w|)}{g(|\nabla w|)}|\nabla w|-1
\Big) 2-2 N\Big\}\\ \di&\geq \ep^3\frac{g(|\nabla w|)}{|\nabla
w|^3} 4\mu^3|x|^2 (4\mu^2|x|^2\delta - K) \geq
\ep^3\frac{g(|\nabla w|)}{|\nabla w|^3} 4\mu^3 {r_2}^2
(4\mu^2{r_2}^2\delta - K)
\end{align*}
where $K=2N$ if $g_0<1$ and $K=2(g_0-1)+2N$ if $g_0>1$. Therefore
if $\mu$ is big enough, depending only on $\delta, g_0, r_2$ and
$N$, we have  $\mathcal{L}w>0$.
\end{proof}


\section{The minimization problem}

In this section we look for minimizers of the  functional
$\mathcal{J}$.
 We begin by discussing the existence of extremals. Next, we prove
 that any minimizer is a subsolution to the equation $\mathcal L u=0$
 and finally, we  prove that $0\le u\le \mbox{sup\,}\varphi_0$.

\begin{teo}\label{existencia}
If $\mathcal{J}(\varphi_0)<\infty$, then there exists a minimizer
of $\mathcal{J}$.
\end{teo}

\begin{proof}
The proof of existence is standard. We write it here for the
reader's convenience and in order to show how the Orlicz spaces and the condition
\eqref{cond} on the function $G$ come into play.

\medskip

Take  a minimizing sequence $(u_n)\subset\K$, then
$\mathcal{J}(u_n)$ is bounded, so $\int_{\Omega}G(|\nabla u_n|)$
and $|\{u_n>0\}|$ are bounded. As  $u_n = \varphi_0$ in
$\partial\Omega$, we have by Lemma~\ref{equi} that $\|\nabla
u_n-\nabla \varphi_0\|_{G}\leq C$ and by Lemma~\ref{poinc} we also
have $\|u_n-\varphi_0\|_{G}\leq C$. Therefore, by Theorem \ref{ref}
 there exists a
subsequence (that we still call $u_n$) and a function $u_0\in
W^{1,G}(\Omega)$ such that
\begin{align*}
&  u_n \rightharpoonup  u_0 \quad \mbox{weakly in }
W^{1,G}(\Omega),
\end{align*}
and by Theorem \ref{imb}
\begin{align*}
&  u_n \rightharpoonup  u_0 \quad \mbox{weakly in }
W^{1,\delta+1}(\Omega),
\end{align*} and by the compactness of the immersions
$W^{1,{\delta+1}}(\Omega)\hookrightarrow L^{\delta+1}(\Omega)$ and
$W^{1,{\delta+1}}(\Omega)\hookrightarrow
L^{\delta+1}(\partial\Omega)$ we have that,
\begin{align*}
& u_n \to u_0 \quad \mbox{a.e. } \Omega.\\
& u_0=\varphi_0\quad\mbox{on}\quad \partial\Omega,\\
\end{align*}
Thus,
\begin{align*}
 |\{&u_0>0\}| \le \liminf_{n\to\infty}|\{u_n>0\}| \quad\mbox{ and}\\
 \int_{\Omega}&G(|\nabla u_0|)\, dx \le \liminf_{n\to\infty}
\int_{\Omega}G(|\nabla u_n|)\, dx.
\end{align*}
In fact,
\begin{equation}\label{G}
\int_\Omega G(|\nabla u_n|)\,dx\ge\int_\Omega G(|\nabla u_0|)\,dx +\int_\Omega g(|\nabla u_0|)\frac{\nabla u_0}
{|\nabla u_0|}\cdot(\nabla u_n-\nabla u_0)\,dx.
\end{equation}

Recall that $\nabla u_n$ converges weakly to $\nabla u_0$ in $L^G$. Now, since  by property $(\widetilde g4)$
$$
\widetilde G\big(g(|\nabla u_0|)\big)\le C G(|\nabla u_0|),
$$

\noindent there holds that $g(|\nabla u_0|)\frac{\nabla u_0}
{|\nabla u_0|}\in L^{\widetilde G}$ so that, by Theorem \ref{ref}
and passing to the limit in \eqref{G} we get
$$
\liminf_{n\to\infty}\int_\Omega G(|\nabla u_n|)\,dx\ge \int_\Omega G(|\nabla u_0|)\,dx.
$$

Hence $u_0\in \K$ and
$$
\mathcal{J}(u_0)\le \liminf_{n\to\infty}\mathcal{J}(u_n) =
\inf_{v\in\K} \mathcal{J}(v).
$$
Therefore, $u_0$ is a minimizer of $\mathcal{J}$ in $\K$.
\end{proof}

{\begin{lema}\label{subsol}  Let $u$ be a
minimizer of
$\mathcal{J}$. Then, $u$ is an $\mathcal{L}$-- subsolution.
\end{lema}}
\begin{proof}
Let $\ep>0$ and $0\leq \xi\in C^{\infty}_0$. Using the minimality
of $u$ and the convexity of $G$ we have
\begin{align*}0 &\leq \frac{1}{\ep} (\mathcal{J}(u-\ep
\xi)-\mathcal{J}(u))\leq \frac{1}{\ep} \int_{\Omega}G(|\nabla
u-\ep \nabla \xi|)-G(|\nabla u|)\, dx\\& \leq \int_{\Omega}-
 g(|\nabla u-\ep \nabla \xi|) \frac{\nabla u-\ep \nabla
\xi}{|\nabla u-\ep \nabla \xi|}\nabla \xi\, dx\end{align*} and if
we take $\ep\rightarrow 0$ we obtain $$0\leq \int_{\Omega}-
g(|\nabla u|) \frac{\nabla u}{|\nabla u|}\nabla \xi\, dx
$$
\end{proof}

{\begin{lema} Let $u$ be a
minimizer of $\mathcal{J}$. Then
$\di 0\leq u\leq \sup_{\Omega}\varphi_0$.
\end{lema}}
\begin{proof}
Let $M=\sup \varphi_0$,  $\ep>0$ and $v=\min(M-u,0)$, then
\begin{align*}0&\leq \frac{1}{\ep}\Big( \mathcal{J}(u+\ep v)-\mathcal{J}(u)\Big)=
\frac{1}{\ep} \Big(\int_{\Omega} G(|\nabla u+\ep \nabla
v|)-G(|\nabla u|)+\lambda\chi_{\{u+\ep v>0\}}
-\lambda\chi_{\{u>0\}}
\, dx \Big)\\
&\leq \frac{1}{\ep} \Big(\int_{\Omega} \Big(G(|\nabla u+\ep
\nabla v|)-G(|\nabla u|)\Big)\, dx\Big)\leq \int_{\Omega}
g(|\nabla u+\ep \nabla v|) \frac{\nabla u+\ep \nabla v}{|\nabla
u+\ep \nabla v|} \nabla v \, dx
\end{align*}
where in the last inequality we are using the convexity of $G$.

Now, takeing $\ep\rightarrow 0$,  using the definition of $v$ and
$(g3)$ we have that,
\begin{align*}0&\leq\int_{\Omega} g(|\nabla u|) \frac{\nabla u}{|\nabla u|}
\nabla v \, dx =-\int_{\{u>M\}} g(|\nabla u|) |\nabla u| \, dx
\leq -\int_{\{u>M\}} G(|\nabla u|)  \, dx\\ &=-\int_{\{u>M\}}
G(|\nabla v|) \, dx,
\end{align*}
therefore $\nabla v=0$ in $\Omega$ and as $v=0$ on $\partial
\Omega$ we have that $v=0$ in $\Omega$ and then $u\leq M$.

To prove that $u\geq 0$ we argue in a similar way. Take
$v=\min(u,0)$, then we have that,
$$0\leq
\frac{1}{\ep} \mathcal{J}(u-\ep v)-\mathcal{J}(u) \leq
-\int_{\Omega} g(|\nabla u-\ep \nabla v|) \frac{\nabla u-\ep
\nabla v}{|\nabla u-\ep \nabla v|} \nabla v \, dx.$$ Therefore
taking $\ep\rightarrow 0$, using the definition of $v$ and $(g3)$
we have that
$$0\geq\int_{\Omega} G(|\nabla v|) \, dx.$$ As in the first
part, we conclude that $u\geq 0$.
\end{proof}

\section{Lipschitz continuity}

In this section we study the regularity of the minimizers of ${\mathcal J}$.
The main result is the local Lipschitz continuity of a minimizer. This result, together with
the rescaling invariance of the minimization problem,
is a key step in the analysis. Once this regularity is proven, a blow up process (passage to the limit in
linear rescalings)  at points of $\partial\{u>0\}$ allows to simplify the analysis by assuming that
$u$ is a plane solution.

As a first step, we prove that minimizers are H\"older continuous.
We use ideas from \cite{DP}, here all the properties of the
function $G$ come into play.

\begin{teo}\label{cotaalpha} For every $0<\alpha<1$,
any minimizer $u$ is in $C^{\alpha}(\Omega)$ and for
$\Omega'\subset\subset\Omega$, $\|u\|_{C^{\alpha}(\Omega')}\leq
C$, where
$C=C(g_0,\delta,\lambda,\|u\|_{\infty},\alpha,\mbox{dist}(\Omega',\partial
\Omega),G(1)).$
\end{teo}
\begin{proof}
We will see that, for every  $0<\alpha<1$ and
$\Omega'\subset\subset\Omega$ there exists $\rho_0$ such that if
$y\in\Omega'$, $0<\rho<\rho_0$ we have that
$$
\frac1{\rho^N}\int_{B_\rho(y)}G(|\nabla u|)\,dx\le C
\rho^{\alpha-1},
$$
for a constant
$C(N,\delta,g_0,\|u\|_{L^\infty(\Omega)},\rho_0,G(1))$.

In fact, let $r>0$ such that, $B_{r}(y)\subset \Omega$. We can
suppose that $y=0$. Then if   $v$ is the solution of
$$
\mathcal{L} v=0\quad\mbox{in }B_r, \qquad v-u\in W_0^{1,G}(B_r),
$$
we have, therefore  by Theorem \ref{desilarga} that
\begin{equation}\label{a1a22}
\int_{B_r} (G(|\nabla u|)-G(|\nabla v|))\, dx \geq C\Big(\int_{
A_2} G(|\nabla u-\nabla v|)  \, dx+\int_{A_1}F(|\nabla u|) |\nabla
u-\nabla v|^2 \, dx\Big),\end{equation} where
$$A_1=\{x\in {B_r}: |\nabla u-\nabla v|\leq 2 |\nabla u|\},\quad
A_2=\{x\in {B_r}: |\nabla u-\nabla v| > 2 |\nabla u|\},$$  and
$C=C( g_0,\delta)$.

On the other hand, by the minimality of $u$, we have
\begin{equation}\label{rn}\int_{B_r} (G(|\nabla u|)-G(|\nabla v|))\, dx
\leq \lambda (|\{v>0\cap {B_r}\}-|\{u>0\cap {B_r}\}|)\leq \lambda
r^N C_N.
\end{equation}
Combining  \eqref{a1a22} and \eqref{rn} we obtain
\begin{align}\label{rn2}&\int_{ A_2} G(|\nabla u-\nabla v|)  \, dx\leq
C\lambda r^N\\\label{rn3} &\int_{A_1}F(|\nabla u|) |\nabla
u-\nabla v|^2 \, dx\leq C\lambda r^N
\end{align}
Let $\ep>0$ and suppose that $r^{\ep}\leq 1/2$. Then, using (g3),
H\"{o}lder's inequality, the definition of $A_1$ and \eqref{rn3}
we obtain,
\begin{equation}\begin{aligned}\label{rn4}\int_{A_1\cap B_{r^{1+\ep}}} G(|\nabla
u-\nabla v&|)\, dx \\&\leq C\Big(\int_{A_1}F(|\nabla u|) |\nabla
u-\nabla v|^2 \, dx\Big)^{1/2} \Big(\int_{ B_{r^{1+\ep}}}
G(|\nabla u|)\,
dx\Big)^{1/2}\\
&\leq C \lambda^{1/2} r^{N/2} \Big(\int_{ B_{r^{1+\ep}}} G(|\nabla
u|)\, dx\Big)^{1/2}.
\end{aligned}\end{equation}
Therefore, by \eqref{rn2} and \eqref{rn4}, we get,

\begin{equation}\label{rn51}\int_{B_{r^{1+\ep}}} G(|\nabla u-\nabla v|)\, dx \leq C \lambda^{1/2}\Big(\lambda^{1/2}r^N+r^{N/2}
\Big(\int_{ B_{r^{1+\ep}}} G(|\nabla u|)\,
dx\Big)^{1/2}\Big).\end{equation}
 On
the other hand
 by property (3) of Lemma \ref{propsol}  we have for every $\beta \in (0,N)$, that there exists a constant
$C=C(\delta,g_0,N,\beta,\|v\|_{L^{\infty}(B_r)})$ such that
\begin{equation}\label{rn8}\int_{B_{r/2}} G(|\nabla v|)\, dx \leq C r^{\beta}.
\end{equation}
By the maximum principle we have,
\begin{equation}\label{cotalinfinito}\|v\|_{L^{\infty}({B_r})}\leq \|v\|_{L^{\infty}(\partial
{B_r})}= \|u\|_{L^{\infty}(\partial
{B_r})}\leq\|u\|_{L^{\infty}({B_r})}\leq
\|v\|_{L^{\infty}({B_r})}\end{equation} where in the last
inequality we are using Lemma \ref{comp}. Then $\|v\|_{L^{\infty}(
{B_r})}=\|u\|_{L^{\infty}( {B_r})}$. This means that the constant
$C$ depends on $\delta,\ g_0,\ N,\ \beta$ and
$\|u\|_{L^{\infty}(B_{r})}$.

\smallskip

By  (G2) we have, $G(|\nabla u|)\leq C(G(|\nabla u-\nabla
v|)+G(|\nabla v|))$. Therefore by \eqref{rn51} and \eqref{rn8},
and for $r\leq 1$ we have,
\begin{align*}\int_{{B_{r^{1+\ep}}}} G(|\nabla u|)\, dx & \leq C
\Big(r^{\beta}(1+\lambda)+\lambda^{1/2}r^{N/2}
\Big({{B_{r^{1+\ep}}}} G(|\nabla u|)\, dx\Big)^{1/2}\Big)\\ &\leq
C \Big(r^{\beta}(1+\lambda)+r^{\beta/2}
(1+\lambda)^{1/2}\Big(\int_{{B_{r^{1+\ep}}}} G(|\nabla u|)\,
dx\Big)^{1/2}\Big).\end{align*} If we call $A=\int_{{B_r^{1+\ep}}}
G(|\nabla u|)\, dx $, we have
\begin{align*}A&\leq
C\Big((1+\lambda) r^{\beta}+(1+\lambda)^{1/2} r^{\beta/2}
A^{1/2}\Big)\leq C\Big((1+\lambda)r^{\beta}
+2(1+\lambda)^{1/2}r^{\beta/2} A^{1/2}\Big)\\&=C \Big(\big(
r^{\beta/2}(1+\lambda)^{1/2}+ A^{1/2}\big)^2-A\Big),\end{align*}
therefore
\begin{align*}
(C+1) A \leq C \Big(r^{\beta/2} (1+\lambda)^{1/2}+ A^{1/2}\Big)^2\\
\Rightarrow (C+1)^{1/2} A^{1/2} \leq C^{1/2}
\Big(r^{\beta/2}(1+\lambda)^{1/2} + A^{1/2}\Big)\\ \Rightarrow
((C+1)^{1/2}-C^{1/2}) A^{1/2} \leq
C^{1/2}r^{\beta/2}(1+\lambda)^{1/2}.
\end{align*}
Thus,  we have the inequality
\begin{equation}\label{clambdafin1}\int_{{B_{r^{1+\ep}}}} G(|\nabla u|)\, dx \leq
((C+1)^{1/2}+C^{1/2})^2 C (1+\lambda)r^{\beta}\end{equation}

Let now,  $0<\alpha<1$, and take  $\ep>0$ such that
$\beta:=(1+\ep)\big(N-(1-\alpha)\big)<N$. Take
$\rho_0=\big(\frac12\big)^{1+1/\ep}$. Then, if $0<\rho<\rho_0$,
taking   $r=\rho^{1/(1+\ep)}$, we have that  $r^\ep<1/2$. And
therefore replacing in \eqref{clambdafin1} we have,
\begin{equation}\label{clambdafin}
\int_{B_\rho}G(|\nabla u|)\le ((C+1)^{1/2}+C^{1/2}) C
(1+\lambda)\rho^{N-(1-\alpha)}
\end{equation}
and by Lemma \ref{calfa} we conclude that for all $0<\alpha<1$,
$u\in C^{\alpha}(B_{\rho})$ for $0<\rho\leq \rho_0$ and
$\|u\|_{C^{\alpha}(B_{\rho})}\leq \overline{C}$ where
$\overline{C}=\overline{C}(N,\alpha,g_0,\delta,\lambda,\rho_0,\|u\|_{L^\infty(\Omega)})$.
\end{proof}

We then have that $u$ is continuous. Therefore, $\{u>0\}$ is open.
We can prove the following property for minimizers.

{\begin{lema}\label{solucion} Let $u$ be a minimizer of
$\mathcal{J}$. Then $u$ is an $\mathcal{L}$--solution in
$\{u>0\}$.
\end{lema}}
\begin{proof}
Let $B\subset \{u>0\}$ and $v$ such that
$$
\begin{cases}
\mathcal{L}v=0 \quad &\mbox{ in } B,\\
v=u \quad &\mbox{ in } B^c.
\end{cases}
$$

By the comparison principle we have that $v\geq u $ in $B$. Thus,
\begin{align*} 0&\geq \int_{\Omega}G(|\nabla u|)-G(|\nabla v|)\,
dx+\lambda |\{u>0\}|-\lambda |\{v>0\}|=\int_{\Omega}G(|\nabla
u|)-G(|\nabla v|)\, dx\\&\geq C\Big(\int_{A_1} F(|\nabla
u|)|\nabla u-\nabla v|^2\, dx+ \int_{A_2} G(|\nabla u-\nabla v|)
\, dx \Big)\end{align*} where we are using Theorem \ref{desilarga}
and $A_1,\ A_2$, and $F$ are as define therein.

Therefore $$ \int_{A_1} F(|\nabla u|)|\nabla u-\nabla v|^2\,
dx=0.$$ Thus, $F(|\nabla u|)|\nabla u-\nabla v|^2=0$ in $A_1$ and,
by the definition of $A_1$, we conclude that $|\nabla u-\nabla
v|=0$ in this set.

On the other hand, we also have
$$\int_{A_2} G(|\nabla u-\nabla v|) \, dx=0$$ so that  $|\nabla u-\nabla
v|=0$ everywhere in $B$.

Hence, as $u=v$ on $\partial B$ we have that $u=v$. Thus,
$\mathcal{L}u=0$ in $B$.

\end{proof}

In order to get the Lipschitz continuity we first  prove the following estimate for minimizers.

\begin{lema}\label{dx}
For all $x\in \Omega$, with $5d(x)< d(x,\partial \Omega)$ we have
$u(x)\leq C d(x)$, where $d(x)=dist(x,\{u=0\})$. The constant $C$ depends
only on $N$ and $\lambda$.
\end{lema}
To prove Lemma \ref{dx} it is enough to prove the following lemma.
In this proof it is essential that the class of functions $G$
satisfying condition \eqref{cond} is closed under the rescaling
$$
G_s(t):=\di\frac {G(st)}{sg(s)}.
$$

{\begin{lema}\label{linf} If $u$ is a minimizer in $B_1$ with
$u(0)=0$, there exists a constant $C$ such that
$\|u\|_{L^{\infty}(B_{1/4})}\leq C$, where $C$ depends only on
$N$,  $\lambda$, $\delta$ and $g_0$.

\end{lema}}
\begin{proof}
Suppose that there exists a sequence $u_k\in\mathcal{K}$ of
minimizers in $B_1(0)$ such that
$$u_k(0)=0\qquad  \mbox{ and } \qquad \di\max_{\overline{B}_{1/4}}u_k(x)>k.$$ Let
$d_k(x)=\di\mbox{dist}(x,\{u_k=0\})$ and
$\mathcal{O}_k=\di\Big\{x\in B_1: d_k(x)\leq
\frac{1-|x|}{3}\Big\}$. Since $u_k(0)=0$ then
$\overline{B}_{1/4}\subset \mathcal{O}_k$, therefore
$$m_k:=\sup_{\mathcal{O}_k}(1-|x|) u_k(x)\geq
 \max_{\overline{B}_{1/4}}(1-|x|) u_k(x)\geq \frac{3}{4} \max_{\overline{B}_{1/4}} u_k(x)
 >\frac{3}{4} k.$$
 For each fix $k$, $u_k$ is bounded, then $(1-|x|)
 u_k(x)\rightarrow 0 \mbox{ when } |x|\rightarrow 1$ which means
 that there exists $x_k\in {\mathcal{O}_k}$ such that $(1-|x_k|)
 u_k(x_k)= \sup_{\mathcal{O}_k}(1-|x|) u_k(x)$, and then
$$u_k(x_k)=\frac{m_k}{1-|x_k|}\geq m_k> \frac{3}{4} k$$
as $x_k\in \mathcal{O}_k$, and $\delta_k:= d_k(x_k)\leq
\frac{1-|x_k|}{3}$. Let $y_k\in \partial\{u_k>0\}\cap B_1$ such
that $|y_k-x_k|=\delta_k$. Then,
$$\begin{array}{ll}\di (1)\ B_{2\delta_k}(y_k)\subset B_1,\\\\ \di \mbox{ since if
} y\in B_{2\delta_k}(y_k) \Rightarrow |y|< 3\delta_k + |x_k|\leq
1,\\\\
\di (2)\ B_{\frac{\delta_k}{2}}(y_k)\subset \mathcal{O}_k, \\\di
\mbox{ since if } y\in B_{\frac{\delta_k}{2}}(y_k) \Rightarrow
|y|\leq \frac{3}{2}\delta_k + |x_k|\leq 1-\frac{3}{2} \delta_k
\Rightarrow d_k(y)\leq \frac{\delta_k}{2}\leq \frac{1-|y|}{3} \ \
\ \mbox{ and }\\\\\di (3) \mbox{ if } z\in
B_{\frac{\delta_k}{2}}(y_k) \Rightarrow 1-|z|\geq
1-|x_k|-|x_k-z|\geq 1-|x_k|-\frac{3}{2} \delta_k\geq
\frac{1-|x_k|}{2}.
\end{array}$$
By (2) we have
$$\max_{\mathcal{O}_k}(1-|x|) u_k(x)\geq \max_{\overline{B_{\frac{\delta_k}{2}}}(y_k)}(1-|x|) u_k(x)\geq
\max_{\overline{B_{\frac{\delta_k}{2}}}(y_k)}\frac{(1-|x_k|)}{2}
u_k(x),$$ where in the last inequality we are using (3). Then,
\begin{equation}\label{uk1}2u_k(x_k)\geq
\max_{\overline{B_{\frac{\delta_k}{2}}}(y_k)}
u_k(x).\end{equation} As $B_{\delta_k}(x_k)\subset \{u_k>0\}$ then
$\mathcal{L}u_k=0$ in $B_{\delta_k}(x_k)$, and by Harnack
inequality in \cite{Li} we have
\begin{equation}\label{uk2}\min_{\overline{B_{\frac{3}{4}\delta_k}}(x_k)} u_k(x)\geq c
u_k(x_k).\end{equation} As
$\overline{B_{\frac{3}{4}\delta_k}}(x_k)\cap
\overline{B_{\frac{\delta_k}{4}}}(y_k)\neq \emptyset$ we have by
\eqref{uk2}
\begin{equation}\label{uk3}
\max_{\overline{B_{\frac{\delta_k}{4}}}(y_k)} u_k(x)\geq c
u_k(x_k).\end{equation} Let
$w_k(x)=\di\frac{u_k(y_k+\frac{\delta_k}{2} x)}{u_k(x_k)}$. Then,
$w_k(0)=0$ and, by \eqref{uk1} and \eqref{uk3} we have,
\begin{align}\label{desiwk}\max_{\overline{B_1}}w_k\leq 2 &\qquad \max_{\overline{B_{1/2}}}
w_k\geq c>0.\end{align} Let now
$$J_k(w)=\int_{B_1} \frac{G(|\nabla w| c_k)}{g(c_k) c_k} \, dx
+\frac{\lambda}{g(c_k) c_k} \int_{B_1} \chi_{\{w>0\}} (x) \, dx$$
where $c_k=\frac{2 u_k(x_k)}{\delta_k}$ so that
$c_k\rightarrow\infty$.

Let us prove, that $w_k$ is a minimizer of $J_k$. In fact, for any
$v\in W^{1,G}(B_1)$ with $v=w_k$ on $\partial B_1$, define $\di
v_k(y)=v\big(\frac{y-y_k}{\delta_k/2}\big) u_k(x_k)$. Thus,
$v_k=u_k$ on $\partial B_{\delta_k/2}(y_k)$. Then,
\begin{align*}J_k(w_k)&=\frac{2^N}{\delta_k^N}\Big(\int_{B_{\frac{\delta_k}{2}}(y_k)}
\frac{G(|\nabla u_k| )}{g(c_k) c_k} \, dy +\frac{\lambda}{g(c_k)
c_k} \int_{B_{\frac{\delta_k}{2}}} \chi_{\{u_k>0\}} (y) \,
dy\Big)\\ &\leq
\frac{2^N}{\delta_k^N}\Big(\int_{B_{\frac{\delta_k}{2}}(y_k)}
\frac{G(|\nabla {v_k}| )}{g(c_k) c_k} \, dy +\frac{\lambda}{g(c_k)
c_k} \int_{B_{\frac{\delta_k}{2}}(y_k)} \chi_{\{{v_k}>0\}} (y) \,
dy\Big)
\\ &=
\int_{B_1} \frac{G(|\nabla {v}| c_k)}{g(c_k) c_k} \, dx
+\frac{\lambda}{g(c_k) c_k} \int_{B_1} \chi_{\{{v}>0\}} (y) \,
dx={J}_k(v).
\end{align*}
Let $g_k(t):=\frac{g(tc_k)}{g(c_k)}$, where the primitive of $g_k$
is $G_k(t)=\frac{G(tc_k)}{g(c_k)c_k}$  and
$\lambda_k=\frac{\lambda}{g(c_k) c_k}\rightarrow 0$. Then,
$$J_k(w)=\int_{B_1} G_k(|\nabla w|) \, dx
+\lambda_k \int_{B_1} \chi_{\{w>0\}} (x) \, dx.$$
 Observe
that for all $k$, $g_k$ satisfies the inequality \eqref{cond},
with the same constants $\delta$ and $g_0$. In fact,
$$\frac{g'_k(t) t}{g_k(t)}= \frac{g'(c_k t) c_k t}{g_k(c_k t)},$$
and then by \eqref{cond} applied to $tc_k$ we have the desirer
inequality.

Let us take $v_k\in W^{1,G}(B_{3/4})$ such that,
\begin{align}\label{eqvk}
\mathcal{L}_k v_k&=0\ \quad \mbox{ in } B_{3/4}\\
v_k&=w_k \quad \mbox{ in } \partial B_{3/4}
\end{align}
where $\mathcal{L}_k$ is the operator associated to $g_k$. By
\eqref{rn51},  \eqref{clambdafin1} (with $\ep=0$ and $r=3/4$) and
the fact that $\lambda_k\rightarrow 0$, we have that
$$\int_{B_{3/4}} G_k(|\nabla w_k-\nabla v_k|)\, dx \leq C
\lambda_k^{1/2},$$ where $C$ depends on $\delta$, $g_0$, $N$ and
$\|w_k\|_{L^{\infty}(B_1)}$ (observe that, since $\lambda_k$ is
bounded for $k$ large then the constants in  \eqref{rn51} and
\eqref{clambdafin1} are independent of $\lambda_k$). We also have,
by \eqref{desiwk}  that $C$ depends only on $\delta$, $g_0$ and
$N$. On the other hand, by (G1) and (g3) we have
$$G_k(t)=\frac{G(tc_k)}{g(c_k)c_k}\geq
\frac{G(c_k)}{(1+g_0)g(c_k)c_k}\min\{t^{g_0+1},t^{\delta+1}\}\geq
\frac{1}{(1+g_0)^2}\min\{t^{g_0+1},t^{\delta+1}\}.$$
  Therefore,
\begin{align*}C\lambda_k^{1/2}&\geq  \int_{B_{3/4}} G_k(|\nabla w_k-\nabla
v_k|)\, dx \\ &\geq \int_{B_{3/4}\cap\{|\nabla w_k-\nabla
v_k|<1\}} \frac{|\nabla w_k-\nabla v_k|^{g_0+1}}{(1+g_0)^2}\,
dx+\int_{B_{3/4}\cap\{|\nabla w_k-\nabla v_k|\geq 1\}}
\frac{|\nabla w_k-\nabla v_k|^{\delta+1}}{(1+g_0)^2}\, dx.
\end{align*}
Hence \begin{equation}\label{deltag0}\begin{aligned}
&A_k:=\int_{B_{3/4}\cap\{|\nabla w_k-\nabla v_k|\geq 1\}} |\nabla
w_k-\nabla v_k|^{\delta+1}\, dx \rightarrow 0 \quad \mbox{ and }\\
&B_k:=\int_{B_{3/4}\cap\{|\nabla w_k-\nabla v_k|< 1\}} |\nabla
w_k-\nabla v_k|^{g_0+1}\, dx \rightarrow 0.
\end{aligned}\end{equation}
By H\"{o}lder inequality and \eqref{deltag0} we have,
$$C_k:=\int_{B_{3/4}\cap\{|\nabla w_k-\nabla v_k|< 1\}} |\nabla
w_k-\nabla v_k|^{\delta+1}\, dx \leq
B_k^{\frac{\delta+1}{g_0+1}}|B_{3/4}|^{\frac{g_0-\delta}{g_0+\delta}}\rightarrow
0,$$ therefore,
\begin{equation}\label{limitedelta}\int_{B_{3/4}} |\nabla w_k-\nabla v_k|^{\delta+1}\,
dx=A_k+C_k\rightarrow 0.\end{equation} As $w_k=v_k$ on $\partial
B_{3/4}$ then
 $p_k=w_k-v_k \in
W^{1,\delta+1}_0(B_{3/4})$
 and by \eqref{limitedelta} we have \begin{equation}\label{ukw1p}p_k
\rightarrow 0 \quad \mbox{ in }
W^{1,\delta+1}_0(B_{3/4}).\end{equation}
  On
the other hand  by Theorem \ref{cotaalpha} we have that,
\begin{equation}\label{cotawk2}\|w_k\|_{C^{\alpha}(B')}\leq
C(\|w_k\|_{L^\infty(B_{3/4})},g_0,\delta,B')\leq
C(g_0,\delta,B')\quad \forall B'\subset\subset
B_{3/4}.\end{equation} (Here again we may
suppose   that the constant $C$ dose not depend on $\lambda_k$,
since $\lambda_k\rightarrow 0$. Also, recall that $\|w_k\|_{L^\infty(B_1)}\le 2$).

 As $v_k$ are solutions of \eqref{eqvk} by Theorem 1.7 in \cite{Li} (see Lemma \ref{calfa}),
we have for $B'\subset\subset B_{3/4}$
\begin{equation}\label{cotavk2}\|v_k\|_{C^{1,\alpha}(B')}\leq
C(N,\delta,g_0,G_k(1),\mbox{dist}(B',\partial
B_{3/4}),\|v_k\|_{L^\infty(B_{3/4})}).\end{equation} But
$G_k(1)=\frac{G(c_k)}{c_kg(c_k)}\le1$ by (g3) and $\|v_k\|_{L^{\infty}(B_{3/4})}\leq
\|w_k\|_{L^{\infty}(\partial B_{3/4})}\leq 2$. Then, this constant
only depends on $N,\delta$ and $g_0$.

Therefore by \eqref{cotawk2} and \eqref{cotavk2} we have that
there exist subsequences, that we call for simplicity $v_k$ and
$w_k$, and functions $w_0$, $v_0$ $\in C^{\alpha}(B')$
for every $B'\subset\subset B_{3/4}$, such that
\begin{align*}
& w_k \rightarrow w_0\quad \mbox{uniformly in } B_{3/4},\\
&v_k \rightarrow v_0 \quad \mbox{uniformly in } B',\\& \nabla v_k
\rightarrow \nabla v_0 \quad \mbox{uniformly in } B'.
\end{align*}
Then,
\begin{align*}
& \nabla w_k \rightarrow \nabla w_0\quad \mbox{ weakly in } L^{\delta+1}(B_{3/4})\\
&p_k=w_k-v_k \rightarrow w_0-v_0 \quad \mbox{ uniformly in } B'.
\end{align*}
But by \eqref{ukw1p} we have $p_k\rightarrow 0$ in
$W^{1,\delta+1}(B')$. Thus, $v_0=w_0$.

Using Harnack's inequality of \cite{Li}, we have that
$$\sup_{B_{1/2}} v_k\leq C \inf_{B_{1/2}} v_k$$ where the constant $C$
depends only on $g_0, \delta, N$. Then, passing to the limit and
using that $v_0=w_0$ we have that
$$\sup_{B_{1/2}} w_0\leq C \inf_{B_{1/2}} w_0.$$
But by \eqref{desiwk},  passing to the limit again, we have that
$\di\sup_{B_{1/2}} w_0>c>0$ and $\di\inf_{B_{1/2}} w_0=0$ since
$w_k(0)=0$ for all $k$, this is a contradiction.
\end{proof}
\begin{proof}[Proof of Lemma \ref{dx}]
Let $x_0\in \{u>0\}$  with $5 d(x_0)<d(x_0,\partial\Omega)$. Take
$\widetilde{u}(x)=\frac{u(y_0+4 d_0 x)}{4 d_0}$, where
$d_0=dist(x_0,\partial\{u>0\})=dist(x_0,y_0)$ with
$y_0\in\partial\{u>0\}$. If we prove that $\widetilde{u}$ is a
minimizer in $B_1(0)$, as $\widetilde{u}(0)=0$ and
$\frac{|x_0-y_0|}{4 d_0}=1/4$,  by Lemma \ref{linf} we have
$$C\geq\widetilde{u}\Big(\frac{x_0-y_0}{4d_0}\Big)=\frac{u(x_0)}{4d_0}
$$ and the result follows.
\\ So, let us prove that $\tilde{u}$ is a
minimizer in
$B_1(0)$. As $5 d(x_0)<d(x_0,\partial \Omega)$ we have, $B_{4
d_0}(y_0)\subset \Omega$. Let $\widetilde{v}\in {W^{1,G}(B_1(0))}$
and $v$ such that $\widetilde{v}(x)=\frac{v(y_0+4 d_0 x)}{4 d_0}$.
Then, changing variables we have,
$$\int_{B_1} G(|\nabla \widetilde{v}|) \, dx = \int_{B_1} G(|\nabla v(y_0+4 d_0 x)|) \, dx =
\int_{B_{4d_0}(y_0)} \frac{G(|\nabla v(y)|)}{d_0^N4^N} \, dy
$$ and $$|\{\widetilde{v}>0\cap B_1\}|= \frac{|\{\widetilde{v}>0\cap
B_{4d_0}(y_0)\}|}{d_0^N4^N}.$$  As $u$ is a minimizer of
$\mathcal{J}$ in $B_{4 d_0}(y_0)$ we have, if
$\widetilde{v}=\widetilde{u}$ on $\partial B_1(0)$,
\begin{align*}\int_{B_{1}(0)} G(|\nabla
\widetilde{u}(x)|)\, dx +\lambda {|\{\widetilde{u}>0\cap
B_{1}(0)\}|}\ &=\int_{B_{4d_0}(y_0)} \frac{G(|\nabla
u(y)|)}{d_0^N4^N}\, dy +\frac{\lambda {|\{{u}>0\cap
B_{4d_0}(y_0)\}|}}{{d_0^N4^N}}\\&
\leq\int_{B_{4d_0}(y_0)}\frac{G(|\nabla v(y)|)}{{d_0^N4^N}}\, dy
+\frac{\lambda  {|\{{v}>0\cap
B_{4d_0}(y_0)\}|}}{{d_0^N4^N}}\\&=\int_{B_{1}(0)} G(|\nabla
{\widetilde v}(x)|) \, dx+\lambda |\{\widetilde{v}>0\cap
B_{1}(0)\}|.
\end{align*}
Therefore, $\widetilde{u}$ is a minimizer of $\mathcal{J}$ in
$B_1(0)$.

\end{proof}

Now we can prove the uniform Lipschtiz continuity of minimizers of
$\mathcal{J}$.

{
\begin{teo}\label{Lip}
Let $u$ be a
minimizer. Then $u$ is locally Lipschitz
continuous in $\Omega$. Moreover, for any connected open subset
$D\subset\subset
 \Omega$ containing free boundary points, the Lipschitz constant of $u$ in $D$ is
estimated by a constant $C$ depending only on $N, g_0, \delta,
dist(D,\partial\Omega)$ and $\lambda$.
\end{teo}}
\begin{proof}
First, take $x$ such that
$d(x)<\frac{1}{5}\mbox{dist}(x,\partial\Omega)$ and
$\widetilde{u}(y)=\frac{1}{d(x)} u(x+d(x)y)$ for $y\in B_1(0)$. By
Lemma \ref{linf} we have $\widetilde{u}(0)\leq C$ in $B_1$, where
$C$ depends only on $N,\lambda,\delta$ and $g_0$. Since $u>0$ in
$B_{d(x)}(x)$, $\mathcal{L}u=0$ in this ball. Thus
$\mathcal{L}\widetilde{u}=0$ in $B_{1}(0)$ By Harnack's inequality
$\widetilde{u}(y)\leq C$ in $B_{1/2}(0)$ where $C$ depends only on
$N,\lambda,\delta$ and $g_0$. Now, by property (2) in Lemma \ref{propsol}, $|\nabla \widetilde{u}(0)|\leq C
\|\widetilde{u}\|_{L^{\infty}(B_{1/2})}\leq C$ where $C$ depends
only on $N,\lambda,\delta$ and $g_0$. Since $\nabla u(x)=\nabla
\widetilde{u}(0)$, the result follows in the case
$d(x)<\frac{1}{5}\mbox{dist}(x,\partial\Omega)$.

Let $r_1$ such that $dist(x,\partial\Omega)\geq r_1>0$ $\forall
x\in D$, take $D'$, satisfying $D\subset\subset D'\subset \subset
\Omega$ given by
$$D'=\{x\in \Omega/ dist(x,D)<r_1/2\}.$$
 If
$d(x)\leq\frac{1}{5}\mbox{dist}(x,\partial\Omega)$ we proved that
$|\nabla u(x)|\leq C$. If
$d(x)>\frac{1}{5}\mbox{dist}(x,\partial\Omega)$, thus $u>0$ in
$B_{\frac{r_1}{5}}(x)$ and $B_{\frac{r_1}{5}}(x)\subset D'$ so
that $|\nabla u(x)|\leq \frac{C}{r_1} \|u\|_{L^{\infty}(D')}.$

To prove the second part of the theorem, consider now any domain
$D$, and $D'$ as in the previous paragraph.  Let us see that
 $\|u\|_{L^{\infty}(D')}$ is bounded by a constant that depends only on
 $N,
D,r_1,\lambda, \delta,$ and  $g_0$ (we argue as in \cite{AC}
Theorem 4.3).
 Let $r_0=\frac{r_1}{5}$, since $D'$ is connected and not
contained in $\{u>0\}\cap\Omega$, there exists $x_0,...,x_k \in
D'$ such that $x_j\in B_{\frac{r_0}{2}}(x_{j-1}) \ j=1,..., k$,
$B_{r_0}(x_{j})\subset\{u>0\}$ $j=0,...,k-1$ and
$B_{r_0}(x_{k})\not \subseteq\{u>0\}$. By Lemma \ref{linf}
$u(x_k)\leq Cr_0$ and by  Harnack's inequality  in \cite{Li} we
have $u(x_{j+1})\geq c u(x_j)$. Inductively we obtain $u(x_0)\leq
C r_0\ \forall x_0\in D'$. Therefore, the supremum of $u$ over
$D'$ can be estimated by a constant depending only on $N,
r_1,\lambda, \delta,$ and  $g_0$.
\end{proof}

Observe that, if we don't use Lemma \ref{dx}, then we obtain that
the Lipschitz constant depends also on
$\|u\|_{L^{\infty}(\Omega)}$ (that is, depends also on the
Dirichlet datum $\varphi_0$).

\section{Nondegeneracy}

In this section we prove the nondegeneracy of a minimizer at the free boundary and
the locally uniform positive density of the sets $\{u>0\}$ and $\{u=0\}$.

\begin{lema}\label{prom2}
Let $\gamma>0$,  $D\subset\subset \Omega$ and $C$ the constant in
Theorem \ref{Lip}. Then, if $C_1>C$, $B_r\subset \Omega$ and $u$
is a
minimizer, there holds that $$ \quad\frac{1}{r}
\Big(\pint_{ B_r} u^{\gamma} \Big)^{1/\gamma}\geq C_1\quad
\mbox{ implies }\quad u>0 \mbox{ in } B_r$$
\end{lema}
\begin{proof}
If $B_r$ contains a free boundary point, as $u$ vanishes at some
point $x_0\in B_r$, and $|\nabla u(x)|\leq C$ in $B_r$, then
$|u(x)-u(x_0)|\leq C r$, that is, $u(x)\leq C r$ in $B_r$ and then
$\frac{1}{r} \Big(\pint_{ B_r} u^{\gamma}\Big)^{1/\gamma} \leq
C$ which is a contradiction.

\end{proof}

\begin{lema}\label{prom1}
For any $\gamma >1$ and for any $0<\kappa<1$ there exists a
constant $c_{\kappa}$ such that, for any
minimizer $u$ and
for every $B_r\subset \Omega$, we have
$$  \quad\frac{1}{r} \Big(\pint_{ B_r} u^{\gamma}\Big)^{1/\gamma} \leq c_{\kappa}
\quad \mbox{ implies }\quad u=0 \mbox{ in } B_{\kappa r},$$ where
$c_{\kappa}$ depends also on $N,\lambda, g_0, \delta$ and $\gamma$.
\end{lema}
\begin{proof}
We may suppose that $r=1$ and that $B_r$ is centered at zero, (if
not, we take the rescaled function $\tilde{u}=\frac{u(x_0+r
x')}{r}$). By Theorem 1.2 in \cite{Li} we have
$$\ep:=\sup_{B_{\sqrt{\kappa}}} u < C \Big(\pint_{B_1} u^{\gamma}
\Big)^{1/\gamma}$$ where $C=C(\kappa,\gamma)$. Now chose $v$
such that
$$v=\begin{cases} C_1 \ep (e^{-\mu |x|^2}-e^{-\mu
{\kappa}^2}) \quad &\mbox{ in } B_{\sqrt{\kappa}}\setminus
B_{\kappa},\\
0 \quad &\mbox{ in } B_{\kappa}. \end{cases}$$
Here the
constants $\mu>0$ and $C_1<0$ are chosen so that
$\mathcal{L}v<0$ in $B_{\sqrt{\kappa}}\setminus B_{\kappa}$ (see
Lemma \ref{exp}) and $v=\ep$ on $\partial B_{\sqrt{\kappa}}$.
Hence, $v\ge u$ on $\partial B_{\sqrt\kappa}$, and therefore if
$$w=
\begin{cases}
\min(u,v)\quad&\mbox{in}\quad B_{\sqrt\kappa},\\
u&\mbox{in}\quad \Omega\setminus B_{\sqrt\kappa},
\end{cases}
$$
$w$ is an admissible function for the minimizing problem. Thus,
using the convexity of $G$, we find that
$$
\begin{aligned}
&\int_{B_\kappa}G(|\nabla u|)\, dx+\lambda|B_\kappa\cap\{u>0\}|\\
&\hskip1cm = \mathcal{J}(u)- \int_{\Omega\setminus
B_\kappa}G(|\nabla u|) \, dx+\lambda|B_\kappa\cap\{u>0\}|
 -\lambda|\Omega\cap\{u>0\}|\\
 &\hskip 1cm \le \mathcal{J}(w)-
 \int_{\Omega\setminus B_\kappa}G(|\nabla u|) \, dx+\lambda |B_\kappa\cap\{u>0\}|
-\lambda|\Omega\cap\{u>0\}|\\
&\hskip1cm=\int_{B_{\sqrt\kappa}\setminus B_\kappa} G(|\nabla
w|)\, dx-\int_{B_{\sqrt\kappa}\setminus B_\kappa}G(|\nabla u|)\,
dx
\\
 &\hskip1cm\le \int_{B_{\sqrt\kappa}\setminus B_\kappa}g(|\nabla w|) \frac{\nabla w}{|\nabla w|}
 (\nabla w-\nabla u)\,dx=- \int_{B_{\sqrt\kappa}\setminus B_\kappa}g(|\nabla w|) \frac{\nabla w}{|\nabla w|}
 \nabla (u-v)^+\,dx\\
 &\hskip1cm =- \int_{(B_{\sqrt\kappa}\setminus B_\kappa)\cap\{u>v\}}g(|\nabla v|)
  \frac{\nabla v}{|\nabla v|}
 \nabla (u-v)^+\,dx
\end{aligned}
$$
and as $v$ is a subsolution we have, $$\begin{aligned}
&\int_{B_\kappa}G(|\nabla u|)\, dx+\lambda|B_\kappa\cap\{u>0\}|\le
-\int_{\partial B_{\kappa}}g(|\nabla v|)
  \frac{\nabla v}{|\nabla v|}u \,\nu
 \,d\mathcal{H}^{N-1}.
\end{aligned}$$
And, as $|\nabla v|\leq C\ep$ we have that
$$\begin{aligned}
&\int_{B_\kappa}G(|\nabla u|)\, dx+\lambda|B_\kappa\cap\{u>0\}|\le
g(C\ep)\int_{\partial B_{\kappa}}u  \,d\mathcal{H}^{N-1}.
\end{aligned}$$
By Sobolev's trace inequality and by $(\widetilde{g}3)$, for
$\widetilde{G}(\alpha)=\lambda$ we have,
$$
\begin{aligned}
&\int_{\partial B_\kappa}u\le C(N,\kappa)\int_{B_\kappa}|\nabla u|+u\, dx\\
&\hskip1cm \le C(N,\kappa)
\Big(\int_{B_\kappa}G\Big(\frac{|\nabla
u|}{\alpha}\Big)+\int_{B_\kappa\cap\{u>0\}}\widetilde{G}(\alpha)+\int_{B_\kappa}u\,
dx \Big)\\&\hskip1cm
 \le C(N,\kappa,\lambda) (1+\ep)\Big(\int_{B_\kappa}G(|\nabla
u|+\lambda |\{u>0\}\cap B_\kappa| \Big)
\end{aligned}
$$
where in the last inequality we are using that $\int_{B_\kappa}u\,
dx\leq \ep |\{u>0\}\cap B_{\kappa}|$. Therefore,
$$\int_{B_\kappa}G(|\nabla u|)\,
dx+\lambda|B_\kappa\cap\{u>0\}|\le g(C\ep) C (1+\ep)
\Big(\int_{B_\kappa}G(|\nabla u|)\,
dx+\lambda|B_\kappa\cap\{u>0\}|\Big).$$
So that, if $\ep$ is
small enough
$$\int_{B_\kappa}G(|\nabla u|)\,
dx+\lambda|B_\kappa\cap\{u>0\}=0.$$
Then, $u=0$ in $B_{\kappa}$ and the
result follows.

\end{proof}
As a corollary we have,

\begin{corol}\label{ucr}
Let $D\subset\subset \Omega$, $x\in D\cap \partial\{u>0\}$. Then
$$\sup_{B_r(x)}u\geq cr,$$ where $c$ is the constant in Lemma
\ref{prom1} corresponding to $\kappa=1/2$ and $\gamma$ fixed.
\end{corol}

\begin{corol}
For any domain $D\subset\subset \Omega$ there exist constants
$c,C$ depending on $N, g_0, \delta, D$ and $\lambda$, such that,
for any
minimizer $u$ and for every $B_r(x)\subset
D\cap\{u>0\}$, touching the free boundary we have
$$ cr \leq u(x) \leq C r$$
\end{corol}
\begin{proof}
It follows by Lemma \ref{dx} and Lemma \ref{prom1}.
\end{proof}
\begin{teo}\label{densprop}
For any domain $D\subset\subset \Omega$ there exists a constant
$c$, with $0<c<1$ depending on $N, g_0, \delta, D$ and $\lambda$,
such that, for any
minimizer $u$ and for every $B_r\subset
\Omega$, centered on the  free boundary we have,
$$ c \leq \frac{|B_r\cap\{u>0\}|}{|B_r|} \leq 1-c$$
\end{teo}
\begin{proof}
First, by Corollary \ref{ucr} we have that there exists $y\in B_r$
such that $u(y)>cr$ and as $u$ is a subsolution we have by Theorem
1.2 in \cite{Li} that
$$\Big(\pint_{B_{\kappa r}} u^{\gamma} \,
dx\Big)^{1/\gamma}\geq C u(y).$$ Therefore
$$\frac{1}{\kappa r}\Big(\pint_{B_{\kappa r}} u^{\gamma} \,
dx\Big)^{1/\gamma}\geq \frac{C}{\kappa}.$$ Now, if $\kappa$ is
small enough, we have
$$\frac{1}{\kappa r}\Big(\pint_{B_{\kappa r}} u^{\gamma} \,
dx\Big)^{1/\gamma}\geq C_1,$$ so that  by Lemma \ref{prom2}, we
have that $u>0$ in $B_{\kappa r},$ where $\kappa=\kappa(C_1,c).$
Thus,
$$\frac{|B_r \cap \{u>0\}|}{|B_r|}\geq \frac{|B_{\kappa
r}|}{|B_r|}= \kappa^N,$$ and $\kappa=\kappa(C_1,c).$

In order to prove the other inequality, we may assume that $r=1$. Let us suppose by contradiction
that, there exists a sequence of minimizers $u_k$ in $B_1$, such
that, $0\in\partial\{u_k>0\}$, with $|\{u_k=0\}\cap
B_{1}|=\ep_k\rightarrow 0$. If we take $v_k\in W^{1,G}(B_{1/2})$
such that,
\begin{align}\label{eqvk2}
\mathcal{L} v_k&=0\ \quad \mbox{ in } B_{1/2}\\
v_k&=u_k \quad \mbox{ in } \partial B_{1/2}
\end{align}

Let $A_1$ and $A_2$  as in  Theorem \ref{desilarga}, for $r=1/2$.
Then we have, by \eqref{rn} that
$$
\begin{aligned}
& \int_{A_2} G(|\nabla u_k-\nabla v_k|)\,
dx\leq C \ep_k\quad \mbox{ and }
\\
& \int_{A_1} F(|\nabla u_k|) |\nabla u_k-\nabla v_k|^2\, dx\leq C
\ep_k,
\end{aligned}
$$
where $C=C(\delta,g_0)$. By \eqref{rn4} (with $\ep=0$ and $r=1/2$)
we have,
$$
\int_{A_1} G(|\nabla u_k-\nabla v_k|)\, dx\leq C \Big(\int_{A_1} F(|\nabla u_k|)
|\nabla u_k-\nabla v_k|^2\, dx\Big)^{1/2} \Big(\int_{A_1}
G(|\nabla u_k|)\Big)^{1/2}.
$$
Therefore, by \eqref{clambdafin1}, there exists $C$ independent of
$k$ such that
$$
\int_{B_{1/2}} G(|\nabla u_k-\nabla v_k|)\,
dx\leq C\ep_k^{1/2}\rightarrow 0.
$$
 As $u_k=v_k$ on $\partial B_{1/2}$, $w_k=u_k-v_k \in W^{1,\delta+1}_0(B_{1/2})$. Thus,
\begin{equation}\label{ukw1p2}
w_k \rightarrow 0 \quad \mbox{  in }
W^{1,\delta+1}_0(B_{1/2}).
\end{equation}

\noindent By Theorem
\ref{cotaalpha} and Theorem 1.7 in \cite{Li}, we have
$$\begin{aligned} \|u_k\|_{C^{\alpha}(B_{1/2})} &\leq
C(N,\delta,g_0,\|u_k\|_{L^{\infty}(B_{1/2})},\alpha) \quad \mbox{
(for } \ep_k \mbox{ small),}\\ \|v_k\|_{C^{1,\alpha}(B')}&\leq
C(N,\delta,g_0,G(1),\|u_k\|_{L^{\infty}(B_{1/2})},\alpha) \quad
\mbox{ (see \eqref{cotalinfinito})}.
\end{aligned}$$
 Therefore,  there exist subsequences, that we call
for simplicity $u_k$ and $v_k$, and functions $v_0\in C^1(B')$,
$u_0$ $\in C(B')$ for all $B'\subset\subset B_{1/2}$ such that
\begin{align*}
& u_k \rightarrow u_0\quad \mbox{ uniformly in } B_{1/2}\\
&v_k \rightarrow v_0 \quad \mbox{ uniformly in } B'\\& \nabla v_k
\rightarrow \nabla v_0 \quad \mbox{ uniformly in } B'
\\& \nabla u_k \rightarrow \nabla u_0\quad \mbox{ weakly in } L^{\delta+1}(B_{1/2})\\
&w_k=u_k-v_k \rightarrow 0 \quad \mbox{ uniformly in } B'.
\end{align*}
 Thus, $v_0=u_0$. By Lemma \ref{prom1} we
have that
$$\Big(\pint_{B_{1/4}} u_k^{\gamma}\Big)^{1/\gamma}\geq C>0.$$
Therefore, passing to the limit, we have
$$\Big(\pint_{B_{1/4}} u_0^{\gamma}\Big)^{1/\gamma}\geq C>0.$$
On the other hand, by Harnack inequality $\sup_{B_{1/4}} v_k\leq
C\inf_{B_{1/4}} v_k$ and again, passing to the limit we have,
$\sup_{B_{1/4}} u_0\leq C\inf_{B_{1/4}} u_0$. As $u_0(0)=0$, then
$u_0\equiv 0$ in $B_{1/4}$, which is a contradiction.

\end{proof}
\begin{remark}\label{rema}
Theorem  \ref{densprop} implies that the free boundary has
Lebesgue measure zero. Moreover, it implies that for every
$D\subset\subset\Omega$, the intersection $\partial\{u>0\}\cap D$
has Hausdorff dimension less than $N$. In fact, to prove these
statements, it is enough to use the left hand side estimate in
Theorem \ref{densprop}. In fact,  this  estimate says that the set
of Lebesgue points of $\chi_{\{u>0\}}$ in $\partial\{u>0\}\cap D$ is empty. On the
other hand almost every point $x_0\in\partial\{u>0\}\cap D$ is a
Lebesgue point, therefore $|\partial\{u>0\}\cap D|=0$.
\end{remark}

\section{The measure $\Lambda=\mathcal{L}u$}
In this section we prove that $\{u>0\}\cap\Omega$ is locally of
finite perimeter. Then, we study the measure
$\Lambda=\mathcal{L}u$ and prove that it is absolutely continuous
with respect to the $\H$ measure on the free boundary. This result
gives rice to a representation theorem for the measure $\Lambda$.
Finally, we prove that almost every point in the free boundary
belongs to the reduced free boundary.

\begin{teo}
For every $\varphi\in C_0^\infty(\Omega)$ such that ${\rm
supp}(\varphi)\subset\{u>0\}$,
\begin{equation}\label{ecuacion}
\int_\Omega g(|\nabla u|)\frac{\nabla u}{|\nabla u|}\nabla
\varphi=0.
\end{equation}
Moreover, the application
$$
\Lambda(\varphi):= -\int_\Omega g(|\nabla u|)\frac{\nabla
u}{|\nabla u|}\nabla \varphi\,dx
$$
from $C_0^\infty(\Omega)$ into $\R$ defines a nonnegative Radon
measure $\Lambda=\mathcal{L}u$ with support on
$\Omega\cap\partial\{u>0\}$.
\end{teo}
\begin{proof}
We know that $u$ is an $\mathcal{L}-$ subsolution, then by the
Riesz Representation Theorem, there exists a nonnegative Radon
measure $\Lambda$, such that $\L u=\Lambda$ . And as
$\mathcal{L}u=0$ in $\{u>0\}$, then for any $\varphi\in
C_0^\infty(\Omega\setminus\partial\{u>0\})$
$$\Lambda(\varphi)=-\int_{\{u>0\}} \nabla \varphi\ g(|\nabla u|) \frac{\nabla u}{|\nabla
u|}\, dx =0,$$ and the result follows.

\end{proof} Now we want to prove that $\Omega\cap
\partial\{u>0\}$, has Hausdorff dimension $N-1$.
First we need the following lemma,

\begin{lema}\label{proplim}
If $u_k$ is a sequence of minimizers in compact subsets of $B_1$,
such that $u_k \rightarrow u_0$ uniformly in $B_1$, then
\begin{enumerate}

\item $\partial\{u_k>0\}\to \partial\{u_0>0\}$ locally in
Hausdorff distance,

\medskip

\item $\chi_{\{u_k>0\}}\to \chi_{\{u_0>0\}}$ in $L^1_{\rm
loc}(\R^N)$,

\medskip

\item If $0\in \partial\{u_k>0\}$, then $0\in
\partial\{u_0>0\}$.
\end{enumerate}
\end{lema}
\begin{proof}
 Here we only have to use
Lemma \ref {prom1} and Theorem \ref{densprop} and the fact that
$u_k\rightarrow u_0$ uniformly in compacts subsets of $B_1$. To
see the complete  proof see pp. 19-20 in \cite{ACF}.
\end{proof}

Now, we prove the following theorem,
\begin{teo}\label{rn-1}
For any domain $D\subset\subset \Omega$ there exist constants $c,
C$,
 depending on $N, g_0, \delta, D$ and $\lambda$, such
that, for any
minimizer $u$ and for every $B_r\subset
\Omega$, centered on the  free boundary we have
$$ c r^{N-1}\leq \int_{B_r}d\Lambda \leq C r^{N-1}$$
\end{teo}

\begin{proof}
Let $\xi\in C_0^\infty(\Omega)$, $\xi\geq 0$. Then,
$$\Lambda(\xi)=-\int_{\{u>0\}} \nabla \xi\, g(|\nabla u|) \frac{\nabla
u}{|\nabla u|}\, dx.$$ Approximating $\chi_{B_r}$ from below by a
sequence  $\{\xi_n\}$ such that $\xi_n=1$ in $B_{r-\frac{1}{n}}$ and
$|\nabla \xi_n|\leq C_N n$ and using that $u$ is Lipschitz we
have that,
$$\Big| \int_{\Omega}\nabla \xi_n\, g(|\nabla u|) \frac{\nabla
u}{|\nabla u|}\, dx \Big|\leq C n \Big|B_{r}\setminus
B_{r-\frac{1}{n}}\Big|\leq C (r^{N-1}+O(1/n)).$$ Then, as
$$\int_{\Omega} \xi_n d\Lambda\to \int_{B_r} d\Lambda,$$
the bound from above holds.

In order to prove the other inequality, we will suppose that
$r=1$. Arguing by contradiction we assume that there exists a
sequence of minimizers $u_k$ in $B_1$, with $0\in
\partial\{u_k>0\}$, and $\Lambda_k=\mathcal{L}u_k$, such that
$\int_{B_1} d\Lambda_k=\ep_k\rightarrow 0$. As the $u_k's$ are
uniformly Lipschitz, we can assume that $u_k\rightarrow u_0$
uniformly in $B_{1/2}$. Let $h_k=g(|\nabla u_k|)\frac{\nabla
u_k}{|\nabla u_k|}$. Then, there exists a subsequence and a
function $h_0$ such that $h_k\rightharpoonup h_0$ $*-$ weakly in
$L^{\infty}(B_{1/2})$. We claim that $h_0=g(|\nabla
u_0|)\frac{\nabla u_0}{|\nabla u_0|}$. In fact, if
$B_{\rho}\subset\subset \{u_0>0\}$ then there exists a subsequence
such that $u_k\rightarrow u_0$ strongly in
$C^{1,\alpha}(B_{\rho})$. So that $h_0=g(|\nabla u_0|)\frac{\nabla
u_0}{|\nabla u_0|}$. If $B_{\rho}\subset \{u_0=0\}$, then by Lemma
\ref{prom1} we have that $u_k=0$ in $B_{\rho\kappa}$ for $k\geq
k_0(\kappa)$. Thus $h_0=0=g(|\nabla u_0|)\frac{\nabla u_0}{|\nabla
u_0|}$ also in this case. Finally $\partial\{u_0>0\}\cap B_{1/2}$
has zero Lebesgue measure. In fact, by (1) in Lemma \ref{proplim},
every point $x_0\in
\partial\{u_0>0\}\cap B_{1/2}$ is a limit point of $x_k\in
\partial\{u_k>0\}\cap B_{1/2}$. Thus, $$\Big(\pint_{B_r(x_0)} u_0^{\gamma}\Big)^{1/\gamma}\geq
cr$$ for any ball $B_r(x_0)\subset B_{1/2}$. Using this fact,
and the Lipschitz continuity we have that
$|B_r(x_0)\cap\{u_0>0\}|\geq c |B_r(x_0)|$ with $c>0$. This implies
that $|\partial\{u_0>0\}\cap B_{1/2}|=0$ (see Remark \ref{rema}).

Therefore, for all $\xi \in C_0^{\infty}(B_{1/2})$, $\xi\geq 0$ we
have
$$\int_{B_{1/2}} g(|\nabla u_0|) \frac{\nabla u_0}{|\nabla u_0|}
\nabla \xi =\lim_{k \rightarrow \infty} \int_{B_{1/2}}  g(|\nabla
u_k|) \frac{\nabla u_k}{|\nabla u_k|} \nabla \xi.
$$ On the other hand,
$$\int_{B_{1/2}}  \xi\, d\Lambda_0 =\lim_{k \rightarrow \infty} \int_{B_{1/2}}
\xi\, d\Lambda_k \leq \|\xi\|_{L^{\infty}(B_{1/2})} \lim_{k
\rightarrow \infty} \ep_k=0.$$ Therefore $\Lambda_0=0$ in
$B_{1/2}$. That is, $\mathcal{L}u_0=0$ in $B_{1/2}$. But $u_0\geq
0$ and  $u_0(0)=0$, so that  by the Harnack inequality we have
$u_0=0$ in $B_{1/2}$.

On the other hand, $0\in \partial\{u_k>0\}$, and by the
nondegeneracy, we have
$$\Big(\int_{B_{1/4}} u_k^{\gamma}\Big)^{1/\gamma}\geq c>0.$$
Thus,
$$\Big(\int_{B_{1/4}} u_0^{\gamma}\Big)^{1/\gamma}\geq c>0$$
which is a contradiction.

\end{proof}
Therefore, we have the following representation theorem

\begin{teo}[Representation Theorem] \label{repteo} Let $u$ be a
minimizer. Then,
\begin{enumerate}
\item $\H( D\cap\partial\{u>0\})<\infty$ for every
$D\subset\subset\Omega$.

\item There exists a Borel function $q_u$ such that
$$
\mathcal{L} u= q_u \,\H \lfloor \partial\{u>0\}.
$$
i.e
$$-\int_\Omega g(|\nabla u|)\frac{\nabla
u}{|\nabla u|}\nabla \varphi\,dx=\int_{\Omega\cap
\partial\{u>0\}} \varphi q_u\,d\mathcal{H}^{N-1} \quad \forall\ \varphi \in \ C_0^\infty(\Omega).
$$
 \item For $D\subset\subset\Omega$ there are constants $0<c\le
C<\infty$ depending on $N,g_0, \delta, \Omega, D$ and $\lambda$
such that for $B_r(x)\subset D$ and $x\in \partial\{u>0\}$,
$$
c\le q_u(x)\le C,\quad c\,r^{N-1}\le
\H(B_r(x)\cap\partial\{u>0\})\le C\,r^{N-1}.
$$

\end{enumerate}
\end{teo}
\begin{proof}
It follows as in Theorem 4.5 in \cite{AC}.
\end{proof}

\begin{remark}\label{perfin}
As $u$ satisfies the conclusions of Theorem \ref{repteo}, the set
$\Omega\cap\{u>0\}$ has finite perimeter locally in $\Omega$ (see
\cite{F} 4.5.11). That is, $\mu_{u}:=-\nabla \chi _{\{u>0\}}$ is a
Borel measure, and the total variation $|\mu_u|$ is a Radon
measure. We define the reduced boundary as in \cite{F}, 4.5.5.
(see also \cite{EG}) by, $\partial_{red}\{u>0\}:=\{x\in
\Omega\cap\partial\{u>0\}/ |\nu_u(x)|=1\}$, where $\nu_u(x)$ is
the unit vector with
\begin{equation}\label{carac}
\int_{B_r(x)}|\chi_{\{u>0\}} -\chi_{\{y/\langle y-x,\nu_u(x)\rangle<0\}}|=o(r^{N})
\end{equation}
for $r\to 0$, if
such a vector exists, and $\nu_u(x)=0$ otherwise. By the results
in \cite{F} Theorem 4.5.6 we have,
$$\mu_u=\nu_u \H \lfloor \partial_{red}\{u>0\}.$$
\end{remark}

\begin{lema}\label{redcasitodo}
$\H(\partial\{u>0\}\setminus\partial_{\rm red}\{u>0\})=0.$
\end{lema}
\begin{proof}
This is a consequence of the density property of Theorem
\ref{densprop} and Theorem 4.5.6 (3) of \cite{F}.
\end{proof}


\section{Asymptotic development and identification of the  function $q_u$}
 In this section we give some properties of blow up sequences of minimizers, we
prove that any limit of a blow up sequence is a minimizer. We
 prove the asymptotic development of minimizers near points
in their reduced free boundary. We finally identify the function
$q_u$ for almost every point in the reduced free boundary.

We first prove some properties of blow up sequences,
\begin{defi}
Let $B_{\rho_k}(x_k)\subset\Omega$ be a sequence of balls with
$\rho_k\to 0$, $x_k\to x_0\in \Omega$ and $u(x_k)=0$. Let
$$
u_k(x):=\frac{1}{\rho_k} u(x_k+\rho_k x).
$$
We call $u_k$ a blow-up sequence with respect to
$B_{\rho_k}(x_k)$.
\end{defi}

Since $u$ is locally Lipschitz continuous, there exists a blow-up
limit $u_0:\R^N\to\R$ such that, for a subsequence,
\begin{align*}
& u_k\to u_0 \quad \mbox{in} \quad C^\alpha_{\rm loc}(\R^N)\quad
\mbox{for every}\quad 0<\alpha<1,\\
& \nabla u_k\to\nabla u_0\quad *-\mbox{weakly  in}\quad
L^\infty_{\rm loc}(\R^N),
\end{align*}
and $u_0$ is Lipschitz in $\RR^N$ with constant $L$.
\begin{lema}\label{propblowup}
If $u$ is a
minimizer then,
\begin{enumerate}

\item $\partial\{u_k>0\}\to \partial\{u_0>0\}$ locally in
Hausdorff distance,

\medskip

\item $\chi_{\{u_k>0\}}\to \chi_{\{u_0>0\}}$ in $L^1_{\rm
loc}(\R^N)$,

\medskip
\item $\nabla u_k\to\nabla u_0$ uniformly in compact subsets of
$\{u_0>0\}$,

\medskip

\item $\nabla u_k\to\nabla u_0$ a.e in $\Omega$,

\medskip

\item If $x_k\in \partial\{u>0\}$, then $0\in
\partial\{u_0>0\}$
\medskip

\item $\mathcal{L}u_0=0$ in
$\{u_0>0\}$.
\end{enumerate}
\end{lema}

\begin{proof}
(1), (2) and (5) follow as in Lemma \ref{proplim}. For the proof
of (3) and (4) see pp. 19-20 in \cite{ACF}. (6) follows by Lemma
\ref{solucion} and by (2) and (3).
\end{proof}
\begin{lema}\label{minblow}
If $u(x_m)=0$, $x_m\rightarrow x_0$ in $\Omega$. Then, any blow up
limit $u_0$ respect to $B_{\rho_m}(x_m)$ is a
minimizer of
$\mathcal{J}$ in any ball.
\end{lema}
\begin{proof}
Let $u_m$, $u_0$ be as is Lemma \ref{propblowup}, $R>0$ and $v$
such that $v-u_0 \in W_0^{1,G}(B_R(0))$. Let $\eta \in
C_0^{\infty}(B_R(0))$, $0\leq \eta\leq 1$ and
$v_m=v+(1-\eta)(u_m-u_0)$ then $v_m=u_m$ in $\partial B_R(0)$.
Therefore
$$\int_{B_R(0)} G(|\nabla u_m|)\, dx + \lambda
\chi_{\{u_m>0\}}\leq \int_{B_R(0)} G(|\nabla v_m|)\, dx + \lambda
\chi_{\{v_m>0\}}.$$ As $|\nabla u_m|\leq C$ and $\nabla
u_m\rightarrow \nabla u_0$ a.e, we have
$$\int_{B_R(0)} G(|\nabla u_m|)\, dx\rightarrow \int_{B_R(0)} G(|\nabla u_0|)\,
dx,$$
$$\int_{B_R(0)} G(|\nabla v_m|)\, dx\rightarrow \int_{B_R(0)} G(|\nabla v|)\,
dx$$ and
$$\chi_{\{v_m>0\}}\leq \chi_{\{v>0\}}+ \chi_{\{\eta<1\}}.$$
Therefore,
$$\int_{B_R(0)} G(|\nabla u_0|)\, dx+\lambda \chi_{\{u_0>0\}}\leq \int_{B_R(0)} G(|\nabla v|)\,
dx, + \lambda \chi_{\{v>0\}}+\lambda \chi_{\{\eta<1\}}.$$ Taking
$\eta$ such that $|\{\eta<1\}\cap B_R(0)|\rightarrow 0$ we have
the desired result.

\end{proof}

Let $\lambda^*$ be such that,
$g(\lambda^*)\lambda^*-G(\lambda^*)=\lambda$. Then we have,

\begin{lema}\label{constant1}
Let $u$ be a
minimizer in $\R^N$ such that $u=\lambda_0\langle
x,\nu_0\rangle^-$ in $B_{r_0}$, with $r_0>0$, $0<\lambda_0<\infty$
and $\nu_0$ a unit vector. Then, $\lambda_0=\lambda^*$.
\end{lema}
\begin{proof}
Let $\tau_{\ep}(x)=x+\ep \eta(x)$ with $\eta\in
C_0^{\infty}(B_{r_0}))$, and let $u_{\ep}(\tau_{\ep}(x))=u(x)$. Then,
$$
0\leq \mathcal{J}(u_{\ep})-\mathcal{J}(u),
$$
$$|B_{r_0}\cap\{u_{\ep}>0\}|=\int_{B_{r_0}\cap\{\langle x,\nu_0\rangle<0\}} |\mbox{det} D\tau_{\ep}|\, dx =
\int_{B_{r_0}\cap\{\langle x,\nu_0\rangle<0\}} (1+\ep \,\mbox{div\,} \eta + o(\ep))\, dx
$$
and
$$
\begin{aligned}
&\int_{B_{r_0}\cap\{u_{\ep}>0\}} G(|\nabla u_{\ep}|) \,
dy \\
&\hskip1cm=\int_{B_{r_0}\cap\{\langle x,\nu_0\rangle<0\}}\Big(G(|\nabla u|) + \ep
\Big(G(|\nabla u|) \mbox{div} \eta - \frac{g(|\nabla u|)}{|\nabla
u|} \nabla u D\eta \nabla u \Big)\Big)\, dx + o(\ep).
\end{aligned}
$$
Therefore, since $u_\ep=u$ in $\RR^N\setminus B_{r_0}$,
$$
0\leq \ep \int_{B_{r_0}\cap\{\langle
x,\nu_0\rangle<0\}}\Big((G(|\nabla u|)+\lambda)  \mbox{div}\eta -
\frac{g(|\nabla u|)}{|\nabla u|} \nabla u D\eta \nabla u \Big)\,
dx + o(\ep).
$$
Thus,
$$ \int_{B_{r_0}\cap\{\langle x,\nu_0\rangle<0\}}\Big((G(|\nabla u|)+\lambda)  \mbox{div}\eta - \frac{g(|\nabla
u|)}{|\nabla u|} \nabla u D\eta \nabla u \Big)\, dx \geq 0.
$$
If we change $\eta$ by $-\eta$ and recall that $\nabla u=-\lambda_0\,\nu_0$ in $\{\langle x,\nu_0\rangle<0\}$ we obtain,
 $$
 \int_{B_{r_0}\cap\{\langle x,\nu_0\rangle<0\}}\Big((G(\lambda_0)+\lambda)  \mbox{div}{\eta} -
g(\lambda_0)\lambda_0 \, \nu_0\,  D{\eta}\, \nu_0 \Big)\, dx =
0
$$
for all $\eta\in C_0^{\infty}(B_{r_0})$.

 Take  $\eta(x)=\phi\displaystyle({|x|})\nu_0$ with
supp\,$\phi\subset (-r_0,r_0)$. Then,
\begin{align*}
&\mbox{div\,}{\eta} (x)=\frac{\phi'(|x|)}{|x|} \langle x,\nu_0\rangle\\
&\nu_0\, D{\eta}\, \nu_0= \nu_{0i}
\frac{\partial{\eta}_j}{\partial x_i} \nu_{0j}=\langle x,\nu_0\rangle
\frac{\phi'(|x|)}{|x|}=\mbox{div\,}{\eta}.
\end{align*}
Hence
\begin{align*}
0&=\int_{\{\langle x,\nu_0\rangle<0\}\cap
B_{r_0}(0)}\big(G(\lambda_0)+\lambda-g(\lambda_0)\lambda_0 \big)
 \mbox{div}{\eta}\, dx \\&= \big(G(\lambda_0)
 +\lambda-g(\lambda_0)\lambda_0 \big)\int_{\{\langle x,\nu_0\rangle=0\}\cap B_{r_0}}
{\eta}\, \nu_0\, d\mathcal{H}^{N-1}(x)\\& =
\big(G(\lambda_0)+\lambda-g(\lambda_0)\lambda_0
\big)\int_{\{\langle x,\nu_0\rangle=0\}\cap B_{r_0}} \phi(|x|)\,
d\mathcal{H}^{N-1}(x)
\end{align*}
for all $\phi \in C_0^{\infty}(-r_0,r_0)$.

Therefore, $g(\lambda_0)\lambda_0-G(\lambda_0)=\lambda$.

\end{proof}

\begin{lema}\label{const}
Let $u\in \K$ be
minimizer. Then, for every  $x_0\in
\Omega\cap\partial\{u>0\}$
\begin{align}
& \limsup_{\stackrel{x\to x_0}{u(x)>0}} |\nabla u(x)| =
\lambda^*.
\end{align}
\end{lema}
\begin{proof}
Let $x_0\in\Omega\cap\partial\{u>0\}$ and let
$$
l:=\limsup_{\stackrel{x\to x_0}{u(x)>0}} |\nabla u(x)|.
$$
Then there exists a sequence $z_k\rightarrow x_0$ such that
$$
u(z_k)>0,\quad \quad |\nabla u(z_k)|\rightarrow l.
$$
Let $y_k$ be the nearest point to $z_k$ on $\Omega \cap
\partial\{u>0\}$ and let $d_k = |z_k-y_k|$. Consider the blow up
sequence with respect to $B_{d_k}(y_k)$ with limit $u_0$, such
that there exists
$$
\nu:=\lim_{k\to\infty} e_k,
$$
 where $e_k=\frac{y_k-z_k}{d_k}$,
and suppose that $\nu = e_N$. Then, by Lemma \ref{propblowup}(1),
$0\in \partial\{u_0>0\}$. By Lemma \ref{propblowup}(2) and by
Lemma \ref{minblow} we have that $u_0$ satisfies Theorem
\ref{densprop}. Then, $B_1(-e_N)\subset \{u_0>0\}$.  By Lemma
\ref{propblowup}(3) we obtain,

  $$|\nabla u_0|\leq l \mbox{ in } \{u_0>0\} \ \mbox{ and }\ |\nabla
 u_0(-e_N)|=l.$$

Then, $0<l<\infty$ and since, by Lemma \ref{propblowup} (6), we
have that $u_0$ is an $\mathcal{L}$ solution in $\{u_0>0\}$ then,
we have that $u$ is locally $C^{1,\alpha}$ there. Thus, there
exists $\mu>0$ such that $|\nabla u_0|>l/2$ in $B_{\mu}(-e_N)$.
Let $e=\frac{\nabla u_0(-e_N)}{|\nabla u_0(-e_N)|}$. Let
$v=\frac{\partial u_0}{\partial e}$, then $v$ satisfies the
uniformly elliptic equation, $D_i(a_{ij}D_j v)=0$ where
$$a_{ij}=\frac{g(|\nabla u_0|)}{|\nabla u_0|}
 \Big[ \Big(\frac{g'(|\nabla u_0|)}{g(|\nabla u_0|)}|\nabla
 u_0|-1\Big) \frac{D_i u_0 D_j u_0}{|\nabla u_0|^2} +\delta_{i
 j}\Big].$$
Then, by the strong maximum principle we have $D_e u_0=l$ in
$B_{\mu}(-e_N)$ so that,   $\nabla u_0=l e$ in $B_{\mu}(-e_N)$. By
continuation we can prove that this is true in $B_1(-e_N)$. Then,
$u_0(x)=l \langle x,e\rangle+C$ in $B_1(-e_N)$. As $u_0(0)=0$  and $u_0>0$
in $B_1(-e_N)$, we have $u_0(x)=l \langle x,e\rangle$ and $e=-e_N$.
Therefore $u_0(x)=-lx_N$ in $B_1(-e_N)$. Using again a
continuation argument we have that $u_0(x)=-lx_N$ in $\{x_N<0\}$.

Now, we want to prove that $u_0=0$ in $\{0<x_N<\ep_0\}$ for some
$\ep_0>0$.

We argue by contradiction. Let
$$
s:=\limsup_{\stackrel{x_N\to 0^+\ x'\in
\mathbb{R}^{N-1}}{u_0(x',x_N)>0}} D_N u_0(x',x_N),
$$
and suppose that $s>0$ ($s<\infty$ since $u_0$ is uniformly
Lipschitz). Let $(z_k,h_k)$ such that, $h_k\rightarrow 0^+$ and
$D_N u_0(z_k,h_k)\rightarrow s$, and take a blow up sequence with
respect to  $B_{h_k}(z_k,0)$ with limit $u_{00}$. Arguing as
before, we have that $u_{00}=sx_N$ for $x_N>0$. On the other hand,
we have $u_{00}=-lx_N$ for $x_N<0$. By Lemma \ref{minblow}
$u_{00}$ is a minimizer, and as all the points of the
form $(x',0)$ belong to the free boundary, we get a contradiction to the positive
density property of the set $\{u_{00}=0\}$(Theorem \ref{densprop}).

Therefore, $s=0$. But this implies that
$u_0(x',x_N)=o(x_N)$ as $x_N\searrow 0^+$. Thus, for all $\ep>0$,
$h_0>0$,
$$
\frac{1}{r} \Big(\pint_{B_r(x_0)} u_0^{\gamma}
\Big)^{1/\gamma} < \ep\ \ \  \mbox{ if } x_0=(y_0,h_0)\ \mbox{and }
r=h_0$$ for $r$ small enough independent of $y_0$. Then, by the
nondegeneracy property, we have that $u_0=0$ in $\{0<x_N<\ep_0\}$.

Now,  by Lemmas \ref{minblow} and \ref{constant1}  we conclude that $l=\lambda^*$, and the result follows.
\end{proof}

\medskip

Now we prove the asymptotic development of minimizers.
 We will  use the following fact.

\begin{remark}\label{lineal} Observe that in  $\{|\nabla u|\geq c\}$, $u$ satisfies a
linear nondivergence uniformly elliptic equation, $Tu=0$ of the
form
\begin{equation}\label{opT}
Tv=b_{ij}(\nabla u) D_{i j}v=0
\end{equation}
where
\begin{equation}\label{bij}
b_{ij}=\delta_{ij}+\Big(\frac{g'(|\nabla u|)|\nabla
 u|}{g(|\nabla u|)}-1\Big) \frac{D_i u D_j u}{|\nabla u|^2},
 \end{equation}
  and the matrix $b_{ij}(\nabla u)$ is
 $\beta$-elliptic in $\{|\nabla
u|>c\}$, where
 $\beta=\max\{\max\{g_0,1\},\max\{1,1/\delta\}\}$.\end{remark}

\begin{teo}\label{blow2}
Let $u$ be a minimizer. Then, at every
$x_0\in\partial_{red}\{u>0\}$, $u$
has the following asymptotic development

 \begin{equation}
 u(x)=\lambda^*\langle x-x_0,\nu(x_0)\rangle^-+o(|x-x_0|).
 \end{equation}
where $\nu(x_0)$ is the outer unit normal to $\partial\{u>0\}$ at $x_0$.
\end{teo}
\begin{proof}
Take $B_{\rho_k}(x_0)$ balls with
$\rho_k\rightarrow 0$ and $u_k$ be a blow up sequence with respect to
these balls with limit $u_0$. Suppose that $\nu_u(x_0)=e_N$.

  First we  prove  that
$$
\begin{cases}
u_0=0 \quad\mbox{ in } \{x_N\geq 0\},\\
u_0>0 \quad\mbox{ in } \{x_N<0\}.
\end{cases}
$$
In fact,
by Lemma \ref{propblowup}, $\chi_{\{u_k>0\}}$ converges to
$\chi_{\{u_0>0\}}$ in $L^1_{loc}$. On the other hand, $\chi_{\{u_k>0\}}$
converges to $\chi_{\{x_N<0\}}$ in $L^1_{loc}$ by
\eqref{carac}. It follows that $u_0=0$ in $\{x_N\geq 0\}$ and
$u_0>0$ a.e in $\{x_N<0\}$.

If $u_0$ were zero somewhere in $\{x_N<0\}$ there should exist a point $\bar x$ in $\{x_N<0\}\cap
\partial\{u_0>0\}$. But, as $u_0$ is a
minimizer, for $0<r<|\bar x_N|$,
$$
\di\frac{|B_r(\bar x)\cap\{u_0=0\}\cap\{x_N<0\}|}{|B_r(\bar x)|}\ge c>0.
$$
Since this is a contradiction we conclude that $u_0>0$ in $\{x_N<0\}$ and therefore
${\mathcal L}u_0=0$ in this set. Since $u_0=0$ on $\{x_N=0\}$, we conclude that
$u_0\in C^{1,\alpha}(\{x_N\le0\})$ (see \cite{Li}). Thus, there exists $0\le\lambda_0<\infty$ such that
$$
u_0(x)=\lambda_0 x_N^-+o(|x|).
$$
By the nondegeneracy of $u$ at every free boundary point (Lemma \ref{prom1}) we deduce that $\lambda_0>0$.

Now, let $u_{00}$ be a blow up  limit of $u_0$. This is, $u_{00}(x)=\lim \frac{u_0(r_n x)}{r_n}$ with $r_n\to0$.
Then, $u_{00}=\lambda_0 x_N^-$. Since $u_{00}$ is
again a
minimizer, Lemma \ref{constant1} gives  that $\lambda_0=\lambda^*$.

Let us see that actually $u_0=\lambda^* x_N^-$. In fact, by applying Lemma \ref{const} we see that
$|\nabla u_0|\le \lambda^*$ and thus, $u_0\le \lambda^* x_N^-$. Since the function $w=\lambda^* x_N^-$ is a solution to
$$
Tw=\sum_{i,j}b_{ij}w_{x_ix_j}=0 \quad\mbox{in }\{x_N<0\}
$$
with $b_{ij}$ as in \eqref{bij} and $u_0$ is a classical solution of the same
equation in a neighborhood of any point where $|\nabla u_0|>0$,  and since $u_0\le w$ in
$\{x_N<0\}$, $u_0=w$ in $\{x_N=0\}$, there holds that either $u_0\equiv w$ or $u_0<w$. In the latter case,
there exists $\delta_0>0$ such that
$$
(w-u_0)(x)\ge -\delta_0\, x_N+o(|x|).
$$
But $(w-u_0)(x)=o(|x|)$. Thus, $u_0\equiv w=\lambda^* x_N^-$.

Finally, since the blow up limit $u_0$ is independent of the blow up sequence $\rho_k$, we deduce that
$$
u(x)=\lambda^*\langle x-x_0,\nu(x_0)\rangle^-+o(|x-x_0|).
$$
\end{proof}

\begin{lema}\label{qu2} For $\mathcal{H}^{N-1}-$ almost every  point
$x_0$ in $\partial_{{red}}\{u>0\}$ there holds that,
$$\int_{B_r(x_0)\cap\partial\{u>0\}}|q_u-q_u(x_0)|
d\mathcal{H}^{N-1}=o(r^{N-1}), \quad \mbox{ as }r\rightarrow 0$$
\end{lema}
\begin{proof}
It follows by Theorem \ref{repteo} (3) that $q_u$ is locally
integrable in $\mathbb{R}^{N-1}$ and therefore almost every point
is a Lebesgue point.
\end{proof}
\begin{lema}
Let $u$ be a
minimizer, then for $\mathcal{H}^{N-1}$ a.e
$x_0\in\partial_{{red}}\{u>0\}$, $$q_u(x_0)=g(\lambda^*).$$

\end{lema}
\begin{proof}
Let $u_0$ be as in Theorem \ref{blow2}.  Now let
 $$\xi(x)=\min\Big(2\big(1-\frac{|x_N|}{2},1\big)\Big)
 \eta(x_1,...,x_{N-1})$$
where $\eta\in C_0^{\infty}(B_r')$, (where $B'_{r}$ is a ball
$(N-1)$ dimensional with radius $r$) and $\eta\geq 0$. Proceeding
as in \cite{AC}, p.121 and using Lemmas \ref{propblowup} and
\ref{qu2}, we get for almost every point
$x_0\in\partial_{red}\{u>0\}$  and
$u_0=\lim_{r\to 0} \frac{u(x_0+rx)}r$ that,
\begin{equation}\label{step3}
-\int_{B_r\cap\{x_N<0\}}g(|\nabla u_0|)\frac{\nabla u_0}{|\nabla
u_0|}\nabla \xi\,dx=q_u(x_0)\int_{ B'_{r}} \xi(x',0)
\,d\mathcal{H}^{N-1} \quad \forall\ \xi \in \ C_0^\infty(B_r),
\end{equation}
where we have assumed that $\nu(x_0)=e_N$.

By Lemma \ref{propblowup}, $u_0=\lambda^* x_N^-$. Substituting in \eqref{step3} we get
$$
g(\lambda^*)\int_{B_r'}\xi(x',0)\,d\mathcal{H}^{N-1}=
q_u(x_0)\int_{ B'_{r}} \xi(x',0)
\,d\mathcal{H}^{N-1} \quad \forall\ \xi \in \ C_0^\infty(B_r).
$$
Thus, $q_u(x_0)=g(\lambda^*)$.
\end{proof}

As a corollary we have
\begin{teo}\label{blow3}
Let $u$ be a
minimizer, then for $\mathcal{H}^{N-1}$ a.e
$x_0\in\partial\{u>0\}$, the following properties hold,

 \begin{equation*}
 q_u(x_0)=g(\lambda^*)
\end{equation*}
and
\begin{equation*}
u(x)=\lambda^*\langle x-x_0,\nu_u(x_0)\rangle^-+o(|x-x_0|)
\end{equation*}

\smallskip

\noindent where $\lambda^*$ is such that,
$g(\lambda^*)\lambda^*-G(\lambda^*)=\lambda$.
\end{teo}
\begin{proof}
The result follows by  Lemma \ref{redcasitodo} and by Theorem
\ref{blow2}.
\end{proof}

\section{Weak solutions}

In this section we introduce the notion of weak solution. The idea, as
in \cite{AC}, is to identify the essential properties that minimizers satisfy
and that may be found in  applications in which minimization does not take
place. For instance, in \cite{MW} we study a singular perturbation problem
for the operator $\mathcal L$ and prove that limits of this singular perturbation
problem are weak solutions in the sense of Definition 8.2.
In the next section, we will prove that weak solutions have smooth free boundaries.
In this way, the regularity results may be applied both to minimizers and to
limits of singular perturbation problems.

With these applications in mind, we introduce two notions of weak solution. Definition
8.1 is similar to the one in \cite{AC} for the case $\mathcal L=\Delta$. On the other
hand, as stated before, Definition 8.2 is more suitable for limits of the singular perturbation
problem.

Since we want to ask as little as possible for a function $u$ to be a weak solution, some
properties already proved for minimizers need a new proof. We keep these proofs as short
as possible by sending the reader to the corresponding proofs for minimizers as soon as possible.

One of the main differences between these two definitions of weak solution is that for
weak solutions according to Definition 8.1 almost every free boundary point is in the reduced
free boundary. Instead, weak solutions according to Definition 8.2 may have an empty
reduced boundary (see, for instance, example 5.8 in \cite{AC}).

In the sequel  $\lambda^*$ will be a fixed positive  constant.

\begin{defi}[Weak solution I]\label{weak}
We call $u$ a weak solution (I), if
\begin{enumerate}
\item $u$ is continuous and non-negative in $\Omega$ and
$\mathcal{L}u=0$ in $\Omega\cap\{u>0\}$. \item For
$D\subset\subset \Omega$ there are constants $0< c_{min}\leq C_{max}$, $\gamma\ge1$, such
that for balls $B_r(x)\subset D$ with $x\in \partial \{u>0\}$
$$
c_{min}\leq \frac{1}{r}\Big(\pint_{B_r(x)} u^{\gamma} dx\, \Big)^{1/\gamma}\leq C_{max}
$$
\item
$$
\mathcal{L} u= g(\lambda^*) \,\H \lfloor
\partial_{red}\{u>0\}.
$$
i.e
$$-\int_\Omega g(|\nabla u|)\frac{\nabla
u}{|\nabla u|}\nabla \varphi\,dx=\int_{\Omega\cap
\partial_{red}{\{u>0\}}} \varphi g(\lambda^*)\,d\mathcal{H}^{N-1} \quad \forall\ \varphi \in \ C_0^\infty(\Omega)
$$
\item
\begin{align*}
& \limsup_{\stackrel{x\to x_0}{u(x)>0}} |\nabla u(x)| \leq
\lambda^*,\qquad \mbox{for every } x_0\in
\Omega\cap\partial\{u>0\}\end{align*}

\end{enumerate}
\end{defi}

\begin{defi}[Weak solution II]\label{weak2}
We call $u$ a weak solution (II), if
\begin{enumerate}
\item $u$ is continuous and non-negative in $\Omega$ and
$\mathcal{L}u=0$ in $\Omega\cap\{u>0\}$. \item For
$D\subset\subset \Omega$ there are constants $0< c_{min}\leq C_{max}$, $\gamma\ge1$, such
that for balls $B_r(x)\subset D$ with $x\in \partial \{u>0\}$
$$
c_{min}\leq \frac{1}{r}\Big(\pint_{B_r(x)} u^{\gamma} dx\, \Big)^{1/\gamma}\leq C_{max}
$$
\item For $\mathcal{H}^{N-1}$ a.e $x_0\in\partial_{{red}}\{u>0\}$,
$u$ has the asymptotic development
$$
u(x)=\lambda^*\langle x-x_0,\nu(x_0)\rangle^-+o(|x-x_0|)
$$
where $\nu(x_0)$ is the unit exterior normal to $\partial\{u>0\}$ at $x_0$ in the measure theoretic sense.

\item
\begin{align*}
& \limsup_{\stackrel{x\to x_0}{u(x)>0}} |\nabla u(x)| \leq
\lambda^*,\qquad \mbox{for every } x_0\in
\Omega\cap\partial\{u>0\}
\end{align*}
\item For any ball $B\subset\{u=0\}$ touching
$\Omega\cap\partial\{u>0\}$ at $x_0$ we have,
$$\limsup_{x\to x_0} \frac{u(x)}{\mbox{dist}(x,B)}\geq  \lambda^*. $$
\end{enumerate}
\end{defi}
\begin{lema}
If $u$ satisfies the hypothesis $(1)$ of Definitions \eqref{weak}
 and \eqref{weak2} then $u$ is in $W^{1,G}_{loc}(\Omega)$ and $\Lambda:=\mathcal{L}u$ is a
nonnegative Radon measure with support in $\Omega\cap\partial\{u>0\}$
(in particular, $u$ is an $\mathcal{L}-$ subsolution in $\Omega$).
\end{lema}
\begin{proof}
Since $\mathcal{L}u=0$ in $\Omega\cap\{u>0\}$, then $u$ is in
$C^{1,\alpha}$ in $\Omega\cap\{u>0\}$. For $s>0$, take
$v=(u-s)^+$. Let $\eta\in
C^{\infty}_0(\Omega)$ with $0\le \eta\leq 1$.  We have,
\begin{align*}
0 &=\int_{\Omega} \frac{g(|\nabla u|)}{|\nabla u|}\nabla u
\nabla (\eta^{g_0+1} v)\, dx\\
&=\int_{\Omega\cap\{u>s\}} \eta^{g_0+1} g(|\nabla u|)|\nabla
u|+ (g_0+1)\int_{\Omega} \eta^{g_0}\, v\,\frac{g(|\nabla u|)}{|\nabla
u|}\nabla u \nabla \eta\, dx.
\end{align*}
Therefore,
\begin{equation}\label{udelta}
\int_{\Omega\cap\{u>s\}} \eta^{g_0+1} g(|\nabla
u|)|\nabla u|\, dx\leq (g_0+1)\int_{\Omega\cap\{u>s\}} g(|\nabla
u|)\, v \, |\eta|^{g_0}|\nabla \eta|,\, dx
\end{equation}
by
\eqref{gmenos3}, ($\widetilde{G}1$) and \eqref{gmenos4} we obtain,
\begin{align*}
\di g(|\nabla u|) |\eta|^{g_0}|v| |\nabla \eta|&\leq
\ep \widetilde{G}(g(|\nabla u|) |\eta|^{g_0})+ C(\ep){G}(|v|
|\nabla \eta|)\\
&\leq C\ep \eta^{g_0+1}\widetilde{G}(g(|\nabla
u|))+ C(\ep){G}(|v| |\nabla \eta|)\\
&\leq C\ep G(|\nabla u|)
\eta^{g_0+1}+C(\ep){G}(|v| |\nabla \eta|).
\end{align*}
Then, using $(g3)$, \eqref{udelta} and choosing $\ep$ small
enough, we have that
$$
\int_{\Omega\cap\{u>s\}}
\eta^{g_0+1} G(|\nabla u|)\, dx\leq C
\int_{\Omega\cap\{u>s\}} {G}(|v| |\nabla \eta|)\, dx\leq C
\int_{\Omega} {G}(|u| |\nabla \eta|)\, dx.
$$
Then, letting $s\to 0$ yields the first assertion.

To prove the second part, take $\xi\in C_0^{\infty}(\Omega)$
nonnegative,  $\ep>0$ and $\di
v=\max\Big(\min\big(1,2-\frac{u}{\ep}\big),0\Big)$. As
$\mathcal{L}u=0$ in $\{u>0\}$, we have that,
\begin{align*}\int_{\Omega}
\frac{g(|\nabla u|)}{|\nabla u|}\nabla u \nabla \xi \,
dx&=\int_{\Omega} \frac{g(|\nabla u|)}{|\nabla
u|}\nabla u \nabla \big(\xi(1-v)\big) \, dx+\int_{\Omega}
\frac{g(|\nabla u|)}{|\nabla u|}\nabla u \nabla (\xi v)\,
dx\\
&=\int_{\Omega} \frac{g(|\nabla u|)}{|\nabla u|}\nabla u
\nabla (\xi v)\, dx= \int_{\Omega\cap\{0<u<2\ep\}}\frac{g(|\nabla u|)}{|\nabla u|}\nabla u
\nabla (\xi v)\, dx\\
&= \int_{\Omega\cap\{\ep<u<2\ep\}}
\frac{g(|\nabla u|)}{|\nabla u|}\nabla u \nabla \Big(\xi
\big(2-\frac{u}{\ep}\big)\Big)\, dx+ \int_{\Omega\cap\{0<u<\ep\}}
\frac{g(|\nabla u|)}{|\nabla u|}\nabla u \nabla \xi\, dx\\
&\leq
2\int_{\Omega\cap\{\ep<u<2\ep\}} g(|\nabla u|) |\nabla \xi|\, dx+
\int_{\Omega\cap\{0<u<\ep\}}g(|\nabla u|)|\nabla \xi |\, dx\\
&\le
2\int_{\Omega\cap\{0<u<2\ep\}}g(|\nabla u|) |\nabla \xi|\, dx,
\end{align*}
which tends to zero when  $\ep\to 0$ yielding  the
desired result.
\end{proof}
Now we will prove as in Theorem \ref{densprop}, the density
property of the set $\{u>0\}$  at free boundary points. It is not true in general, for weak solutions
satisfying only properties (1) and (2) of Definitions \ref{weak}
or \ref{weak2}  that the set $\{u=0\}$ has positive
density at $\H-$ almost every free boundary point (see examples in \cite{AC}).

\begin{teo}\label{densprop2}
For any domain $D\subset\subset \Omega$ there exists a constant
$c$, with $0<c<1$ depending on $N,\gamma, g_0, \delta, D$,  $c_{min}$ and
$C_{max}$, such that, for any function $u$ satisfying (1)
and (2) of Definitions \ref{weak} and \ref{weak2} and for every
$B_r\subset D$, centered at the free boundary we have
$$\frac{|B_r\cap\{u>0\}|}{|B_r|}\geq c$$
\end{teo}
\begin{proof}
The proof follows as in Theorem \ref{densprop}, the only
difference here is that,  instead of using Lemma \ref{prom2}
and \ref{prom1}, we use property (2) of Definitions \ref{weak} and
\ref{weak2}.
\end{proof}
\begin{remark}
Now, by Remark \ref{rema} we have that the free boundary has
Lebesgue measure zero. Moreover,  for every
$D\subset\subset\Omega$, the intersection $\partial\{u>0\}\cap D$
has Hausdorff dimension less than $N$.
\end{remark}

\begin{lema}
If $u$ satisfies  hypothesis $(1)$ and $(2)$ of Definitions
\eqref{weak} and \eqref{weak2} then\begin{enumerate} \item
 $u$ is Lipschitz and for any domain
$D\subset\subset\Omega$, the Lipschitz constant depends only on
$N,\gamma, g_0,\delta, dist(D,\partial\Omega)$ and $C_{max}$, provided $D$
contains a free boundary point. \item For any domain
$D\subset\subset \Omega$ there exist constants $\,c, C$
 depending on $N,\gamma, g_0, \delta, D$,  $c_{min}$ and $C_{max}$, such
that,  for every $B_r\subset D$
centered at the  free boundary we have
$$ c r^{N-1}\leq \int_{B_r}d\Lambda \leq C r^{N-1}.$$
\end{enumerate}
\end{lema}
\begin{proof}
The proof of (1) is similar to the one in Theorem \ref{Lip}. The
only change that we have to  make here  is the following, instead
of using Lemma \ref{dx} we have to use property $(2)$ of
Definitions \ref{weak} and \ref{weak2}. We give the proof for the
readers convenience.

Let $d(x)=\mbox{dist}(x,\Omega\cap\partial\{u>0\})$.  First,
take $x$ such that
$d(x)<\frac{1}{5}\mbox{dist}(x,\partial\Omega)$. Let $y\in
\partial\{u>0\}\cap\partial B_{d(x)}(x)$. As $u>0$ in $B_{d(x)}(x)$, $\mathcal{L}u=0$
in that ball and $u$ is an $\mathcal{L}-$ subsolution in
$B_{3d(x)}(y)$. By using the gradient estimates and Harnack's
inequality of \cite{Li} (see Lemma \ref{propsol}) and   property $(2)$ of Definitions
\ref{weak} and \ref{weak2} we have,
$$|\nabla u(x)|\leq C\frac{1}{d(x)} \sup_{B_{d(x)}(x)} u \leq
C\frac{1}{d(x)} \sup_{B_{2d(x)}(y)} u\leq C\frac{1}{d(x)}
\Big(\pint_{B_{3d(x)}(y)} u^{\gamma} dx\, \Big)^{1/\gamma}\leq
CC_{max}.$$

So, the result follows in the case
$d(x)<\frac{1}{5}\mbox{dist}(x,\partial\Omega)$.

Let $r_1$ such that $dist(x,\partial\Omega)\geq r_1>0$ $\forall
x\in D$, take $D'$, satisfying $D\subset\subset D'\subset \subset
\Omega$ given by
$$D'=\{x\in \Omega/ dist(x,D)<r_1/2\}.$$

Let $x\in D$. If
$d(x)\leq\frac{1}{5}\mbox{dist}(x,\partial\Omega)$ we have proved that
$|\nabla u(x)|\leq C$.

 If
$d(x)>\frac{1}{5}\mbox{dist}(x,\partial\Omega)$,  $u>0$ in
$B_{\frac{r_1}{5}}(x)$ and $B_{\frac{r_1}{5}}(x)\subset D'$ so
that $|\nabla u(x)|\leq \frac{C}{r_1} \|u\|_{L^{\infty}(D')}.$

To prove the second part of (1), consider now a connected domain $D$ that contains a free
boundary point and let
$D'$ as in the previous paragraph.  Let us see that
 $\|u\|_{L^{\infty}(D')}$ is bounded by a constant that depends only on
 $N,\gamma,
D,r_1,\lambda, \delta,$ and  $g_0$.
 Let $r_0=\frac{r_1}{4}$ and $x_0\in D$. Since $D'$ is connected and not
contained in $\{u>0\}\cap\Omega$, there exists $x_1,...,x_k \in
D'$ such that $x_j\in B_{\frac{r_0}{2}}(x_{j-1}) \ j=1,..., k$,
$B_{r_0}(x_{j})\subset\{u>0\}$ $j=0,...,k-1$ and
$B_{r_0}(x_{k})\not \subseteq\{u>0\}$. Let $y_0\in
\partial\{u>0\}\cap B_{r_0}(x_k)$. As $u$ is an $\mathcal{L}-$ subsolution,
by Theorem 1.2 in  \cite{Li} there exists $C$ depending on
$N, \gamma,\delta, g_0$ such that,
$$u(x_k)\leq C \Big(\pint_{B_{2r_0}(y_0)} u^{\gamma} dx\, \Big)^{1/\gamma}\leq
CC_{max}r_0,$$ where in the last inequality we have used property
$(2)$ of Definitions \ref{weak} and \ref{weak2}. By Harnack's
inequality in \cite{Li} we have $u(x_{j+1})\geq c u(x_j)$.
Inductively we obtain $u(x_0)\leq C r_0\ \forall x_0\in D'$.
Therefore, the supremum of $u$ over $D'$ can be estimated by a
constant depending only on $N, \gamma, r_1,\lambda, \delta,$ and  $g_0$.

In order to prove (2) we use that  Lemma \ref{proplim} holds if
$u_k$ is a sequence of functions satisfying properties (1) and (2)
of Definitions \ref{weak} and \ref{weak2} with the same constants
$c_{min}$ and $c_{max}$. Then, the rest of the proof follows as in
Theorem \ref{rn-1}.
\end{proof}

\begin{remark}\label{repweak}Now, we are under  the  conditions used in the proof of Theorem \ref{repteo}
and therefore this result applies to functions  $u$ satisfying
properties $(1)$ and $(2)$ of Definition \ref{weak} and
\ref{weak2}. This is, $\Omega\cap\partial\{u>0\}$ has finite perimeter and
there exists a Borel function $q_u$ defined on $\Omega\cap\partial\{u>0\}$ such that
${\mathcal L}u=q_u\H \lfloor \partial\{u>0\}$.

As $u$ satisfies the conclusions of Theorem
\ref{repteo} then Remark \ref{perfin} also holds. We also have
that any blow up  sequence satisfies the properties  of Lemma
\ref{propblowup}.
\end{remark}
Moreover, we have the following result that  holds at points $x_0\in
\partial_{red}\{u>0\}$ that are Lebesgue points of the function $q_u$ and are such that
\begin{equation}\label{limsupprop}
\limsup_{r\to 0}\frac{\H(\partial\{u>0\}\cap B(x_0,r))}{\H(B'(x_0,r))}\leq
1.
\end{equation}
(Here $B'(x_0,r)=\{x'\in \RR^{N-1}\,/\, |x'|<r\}$).

\smallskip

Recall that
$\H -\,a.e.$ point in $\partial_{red}\{u>0\}$ satisfies \eqref{limsupprop} (see Theorem 3.1.21 in \cite{F}).

\begin{lema}\label{glambda*}
If $u$ is a function satisfying properties (1), (2) and (3) of Definition
\ref{weak} or \ref{weak2} we have that $q_u(x_0)=g(\lambda^*)$ for
$\mathcal{H}^{N-1}$ a.e  $x_0\in
\partial_{red} \{u>0\}$.
\end{lema}
\begin{proof}
Clearly, we only have to prove the statement for weak solutions (II).

If $u$ satisfies (3) of Definition \ref{weak2}, take $x_0\in
\partial_{red}\{u>0\}$ such that
$$
u(x)=\lambda^*\langle x-x_0,\nu(x_0)\rangle^-+o(|x-x_0|).
$$
Take
$\rho_k\to 0$ and $u_k(x)=\frac{1}{\rho_k}u(x_0+\rho_k x).$  If
$\xi\in C_0^{\infty}(\Omega)$ we have
$$-\int_{\{ u>0\}} g(|\nabla u|) \frac{\nabla u}{|\nabla
u|}\nabla \xi\, dx=\int_{\partial \{u>0\}} q_u(x)  \xi d\H,$$ and
if we replace $\xi$ by $\xi_k(x)=\rho_k\xi(\frac{x-x_0}{\rho_k} )$
with $\xi\in C_0^{\infty}(B_R)$, $k\geq k_0$ and we change
variables we obtain,
$$
-\int_{\{ u_k>0\}} g(|\nabla u_k|) \frac{\nabla u_k}{|\nabla
u_k|}\nabla \xi\, dx=\int_{\partial \{u_k>0\}} q_u(x_0+\rho_k x)
\xi d\H.
$$

Now, recall that for a subsequence,
 $\chi_{\{u_k>0\}}\to
\chi_{\{x_N<0\}}$ in $L^1_{\rm loc}(\R^N)$ and
$ g(|\nabla u_k|) \frac{\nabla u_k}{|\nabla
u_k|}\rightharpoonup g(|\nabla u_0|) \frac{\nabla u_0}{|\nabla
u_0|}$ $*-$ weakly in $L^\infty_{loc}(\RR^N)$. Thus,
$$
\int_{\{ u_k>0\}} g(|\nabla u_k|) \frac{\nabla u_k}{|\nabla
u_k|}\nabla \xi\, dx \to
\int_{\{x_N<0\}}g(|\nabla u_0|) \frac{\nabla u_0}{|\nabla
u_0|}\nabla \xi\, dx
$$

On the other hand, $\partial\{u_k>0\}\to \{x_N=0\}$ locally in
Hausdorff distance.
Then, if $x_0$ is a Lebesgue point of $q_u$ satisfying \eqref{limsupprop},
\begin{equation}\label{convergencia}
\int_{\partial \{u_k>0\}} q_u(x_0+\rho_k x) \xi\, d\H\rightarrow
q_u(x_0)\int_{\{x_N=0\}} \xi \,d\H.
\end{equation}

As, $\nabla u_0=-\lambda^* e_N \chi_{\{x_N<0\}}$, we deduce that for almost every point
$x_0 \in \partial_{red}\{u>0\}$,
$q_u(x_0)=g(\lambda^*)$.
\end{proof}

Now we prove the asymptotic development for weak solutions
satisfying Definition \ref{weak}.

\begin{lema}
\label{blow4} If  $u$ satisfies $(1), (2)$, $(3)$ and $(4)$ of
Definition \ref{weak}, then for $x_0\in\partial_{{red}}\{u>0\}$
satisfying  \eqref{limsupprop}, $u$ has the following asymptotic
development
\begin{equation}\label{asym}
u(x)=\lambda^* \langle x-x_0,\nu(x_0)\rangle^-+o(|x-x_0|)
\end{equation}
where  $\nu(x_0)$ is the unit outer normal
to the free boundary at $x_0$.
  \end{lema}
\begin{proof}
Let $x_0\in\partial_{red}\{u>0\}$ and let $\rho_k\to0$. Let
$u_k(x)=\frac1{\rho_k}u(x_0+\rho_k x)$ be a blow up sequence
(observe that $u_k$ is again a weak solution in the rescaled domain). Assume that $u_k\to
u_0$ uniformly on compact subsets of $\RR^N$. Also assume that
$\nu(x_0)=e_N$. As in the proof of Theorem \ref{blow2} we deduce
that
$$
\begin{aligned}
&u_0\ge0\quad\mbox{in }\{x_N<0\}\\
&u_0=0\quad\mbox{in }\{x_N\geq0\}.
\end{aligned}
$$
Let us see that $u_0>0$ in $\{x_N<0\}$. To this end, let $D\subset\subset\{x_N<0\}$ and
let $\xi\in C_0^\infty(D)$. For $k$ large enough,
\begin{equation}\label{uk}
-\int_{\{u_k>0\}} g(|\nabla u_k|)\di\frac{\nabla u_k}{|\nabla u_k|}\nabla\xi\,dx=
\int_{\partial_{red}\{u_k>0\}}g(\lambda^*)\xi(x)\, d{\mathcal H}^{N-1}.
\end{equation}

As in \cite{AC}, p. 121, we have that for every $x_0\in\partial_{red}\{u>0\}$ satisfying \eqref{limsupprop},
$$
{\mathcal H}^{N-1}(\partial\{u_k>0\}\cap D)\to0\quad\mbox{as }k\to\infty.
$$
Thus,  the right hand side of
\eqref{uk} goes to zero as $k\to\infty$. Since the left hand side goes to
$$
-\int g(|\nabla u_0|)\di\frac{\nabla u_0}{|\nabla u_0|}\nabla\xi\,dx
$$
we deduce that ${\mathcal L}u_0=0$ in $\{x_N<0\}$. Thus, $u_0>0$ in
$\{x_N<0\}$.

As in Theorem \ref{blow2} we have that there exists $0<\lambda_0<\infty$ such
that
$$
u_0(x)=\lambda_0 x_N^-+o(|x|).
$$
By property (2) of Lemma \ref{propblowup} we have  that
$$
\chi_{\{u_k>0\}}\to \chi_{\{x_N<0\}}\quad\mbox{in }L^1_{loc}(\RR^N)\quad\mbox{as }k\to\infty.
$$

Let now $\xi\in C_0^\infty(\RR^N)$ in \eqref{uk}. Passing to the
limit as $k\to\infty$ and using Lemma \ref{propblowup} (1) we get,
$$
-\int_{\{x_N<0\}} g(|\nabla u_0|)\di\frac{\nabla u_0}{|\nabla u_0|}\nabla\xi\,dx=
\int_{\{x_N=0\}}g(\lambda^*)\xi(x)\, d{\mathcal H}^{N-1}.
$$

Replacing $\xi$ by $r\xi(x/ r)$ with $r\to0$, using the fact that
$\frac 1r u_0(rx)\to\lambda_0 x_N^-$ uniformly on compact sets of
$\RR^N$, changing variables and passing to the limit we get
$$
g(\lambda_0)\int_{\{x_N<0\}} \xi_N\,dx= g(\lambda^*)
\int_{\{x_N=0\}}\xi(x)\, d{\mathcal H}^{N-1}.
$$
Thus, $\lambda_0=\lambda^*$.

At this point we proceed as in Theorem \ref{blow2} to deduce that
actually $u_0(x)=\lambda^* x_N^-$ (observe that here we are using
property (4) of Definition \ref{weak}). As the blow up limit $u_0$
is independent of the blow up sequence $\rho_k$ we conclude that
$u$ has the asymptotic development \eqref{asym}.
\end{proof}
Now we prove the property that we mentioned in the introduction to this
section. The following lemma only holds  for weak solutions
satisfying Definition \ref{weak}.

\begin{lema}\label{prop2weak}
If $u$ satisfies $(1)$, $(2)$ and $(3)$ of Definition \ref{weak},
\begin{enumerate}
\item $\H(\partial\{u>0\}\setminus \partial_{red}\{u>0\})=0$ \item
$|D\cap\{u=0\}|>0$ for every open set $D\subset\Omega$ containing
a point of $\{u=0\}$. \item For any ball $B$ in $\{u=0\}$ touching
$\Omega\cap\partial\{u>0\}$ at $x_0$, there holds that,
\begin{equation}\label{limsup22}
\limsup_{x\to x_0} \frac{u(x)}{\mbox{dist}(x,B)}\geq  \lambda^*
\end{equation}
 \end{enumerate}
\end{lema}
\begin{proof} By \cite{F}, 4.5.6 (3) we have,
\begin{equation}\label{nored}|\mu_u|(B_r(x_0))=o(r^{N-1})\ \ \mbox{ for } r\to 0\end{equation} for
$\H$ almost all points $x_0\in \partial\{u>0\}\setminus
\partial_{red}\{u>0\}$ (Recall that $\mu_u=-\nabla\chi_{\{u>0\}}$) . Assume there exists
$x_0 \in \partial\{u>0\}\setminus
\partial_{red}\{u>0\}$ satisfying \eqref{nored}. Therefore, if
$u_0$ is a blow up limit with respect to balls $B_{\rho_k}(x_0)$,
we obtain for $\xi\in C_0^{\infty}(B_1)$ that,
\begin{align*}-\int_{\R^N}g(|\nabla u_0|)\frac{\nabla u_0}{|\nabla
u_0|}\nabla \xi\,dx \leftarrow &-\int_{\R^N}g(|\nabla
u_k|)\frac{\nabla u_k}{|\nabla u_k|}\nabla \xi\,dx
\\&=\rho_k^{1-N}g(\lambda^*)\int_{ \partial_{red}\{u>0\}\cap B_{\rho_k}(x_0)}
\xi\Big(\frac{y-x_0}{\rho_k}\Big) \,d\mathcal{H}^{N-1}\\&
=\rho_k^{1-N}g(\lambda^*)\int_{ B_{\rho_k}(x_0)}
\xi\Big(\frac{y-x_0}{\rho_k}\Big) \,d|\mu_u|(x)\\&\leq C
\rho_k^{1-N} |\mu_u|(B_{\rho_k}(x_0))\to 0,
\end{align*}
therefore $\mathcal{L}u_0=0$. Since $u_0(0)=0$, we must have
$u_0=0$, but this contradicts the   nondegeneracy property (2) of
the Definition \ref{weak}. Therefore $(1)$ holds.

To prove $(2)$, suppose that $\chi_{\{u>0\}}=1$ almost everywhere
in $D$, hence the reduced boundary must be outside of $D$. Then by
Definition \ref{weak} (3) the function $\mathcal{L}u=0$ in $D$,
and therefore $u$ is positive. Hence $D\cap\{u=0\}=\emptyset$.

In order to prove (3), Let $l$ be the finite limit on the left of
\eqref{limsup22}, and $y_k \to x_0$ with $u(y_k)>0$ and
$$\frac{u(y_k)}{d_k}\to l, \quad d_k=\mbox{dist}(y_k,B).$$
Consider the blow up sequence $u_k$ with respect to
$B_{d_k}(x_k)$, where $x_k\in\partial B$ are points with
$|x_k-y_k|=d_k$, and choose a subsequence with blow up limit
$u_0$, such that
$$e:=\lim_{k\to\infty} \frac{x_k-y_k}{d_k}$$ exists. Then by
construction, since $l>0$ by nondegenaracy, $u_0(-e)=l$, and
$u_0(x)\leq -l\langle x,e\rangle$ for $x \cdot e\leq 0$,
$u_0(x)=0$ for $x\cdot e\geq 0.$ Both, $u_0$ and $l\langle x,
e\rangle^-$ are $\L$ solutions in $\{u_0>0\}$, and coincide in
$-e$. Since $l>0$, and $|\nabla u_0|>l/2$ in a neighborhood of
$-e$, we have that $\L$ is uniformly elliptic there. Then we can
apply the strong maximum principle to conclude that they must
coincide in that neighborhood of $-e$. By a continuation argument,
we have that $u_0=l\langle x, e\rangle^-$.

By the Representation Theorem, $\forall\ \varphi \in \
C_0^\infty(B_1)$, $\varphi\geq 0$
\begin{equation}
\label{qukk1}\begin{aligned} \int_{
\partial\{u_k>0\}} \varphi q_{u_k}\,d\mathcal{H}^{N-1}=-\int_{\R^N} g(|\nabla u_k|)\frac{\nabla
u_k}{|\nabla u_k|}\nabla \varphi\,dx \rightarrow
&-\int_{\mathbb{R}^N} g(|\nabla u_0|)\frac{\nabla u_0}{|\nabla
u_0|}\nabla \varphi\,dx\\&=g(l)\int_{\{\langle x, e\rangle=0\}}
\varphi\,d\H
\end{aligned}\end{equation}
 and
\begin{equation}\label{qukk2}\begin{aligned} \int_{
\partial\{u_k>0\}} \varphi \,d\mathcal{H}^{N-1}&\geq
\int_{
\partial_{red}\{u_k>0\}} \varphi \langle e. \nu_{u_k}\rangle\,d\mathcal{H}^{N-1}\\&=
\int \varphi e. d\mu_{u_k}=\int_{\{u_k>0\}}\partial_e \varphi\,
dx\rightarrow \int_{\{\langle x, e\rangle<0\}}\partial_e \varphi\,
dx\\&=\int_{\{\langle x, e\rangle=0\}} \varphi d\H.
\end{aligned}\end{equation}

Therefore, for weak solutions of type I and II we have,
$$g(l)\geq \liminf_{x\to x_0} q_u(x).$$

Now, if  $u$ is a weak solution of type I we have, that
$q_u(x)=g(\lambda^*)$ for $\H-$ a.e $x\in \Omega \cap
\partial\{u>0\}$. Thus, $g(l)\geq g(\lambda^*)$ and $l\geq
\lambda^*$.

\end{proof}
We then conclude,
  \begin{teo}\label{toda}
If  $u$ satisfies $(1), (2)$ $(3)$ and $(4)$ of Definition
\ref{weak}, then for $\mathcal{H}^{N-1}$ a.e
$x_0\in\partial\{u>0\}$, $u$ has the asymptotic development \eqref{asym}
\end{teo}
\begin{proof}
It follows by Remark \ref{repweak} and Lemmas \ref{blow4} and \ref{prop2weak}.
\end{proof}

\begin{remark}\label{diferencia} Now we have that with the additional
hypothesis $(4)$, weak solutions (I) satisfy  the same  properties
that we proved in the previous section for minimizers (with the
only difference that in (4) we have a less than or equal instead
of an equal). The extra hypothesis  $(5)$, in the definition of
weak solution (II) (which  always holds, by Lemma \ref{prop2weak},
for weak solutions (I))  is  used in key steps of the proof of the regularity of
the free boundary. On the other hand, observe that minimizers have the asymptotic
development \eqref{asym} at every point in their reduced free
boundary, but we only proved that this development holds at almost
every point of $\partial_{red}\{u>0\}$ when $u$ is a weak
solution.

\end{remark}

\section{Regularity of the free boundary}
In this section we prove the regularity of the free boundary of a
weak solution $u$ in a neighborhood of every ``flat'' free
boundary point. In particular, we prove the regularity in a
neighborhood  of every point in $\partial_{red}\{u>0\}$ where $u$
has the asymptotic development \eqref{asym}. Then, if $u$ is a
minimizer, $\partial_{red}\{u>0\}$ is smooth and the remainder of
the free boundary has $\H-$  measure zero.

We will recall some definitions and we will point
out the only significant differences with the proofs in \cite{DP}
for the case $G(t)=t^p$. The rest of the proof of the regularity
then follows as sections 6, 7, 8 and 9 of \cite{DP}.
The main differences with \cite{DP} come from the fact that we
don't assume the locally uniform positive density of the set $\{u\equiv0\}$ at
the free boundary. This is a property satisfied by minimizers that
is not know to hold, in principle, for weak solutions that appear in a different
context. This uniform density property implies, in particular, that $\H-$ almost every
point in the free boundary belongs to the reduced free boundary and
this is a very strong assumption that we don't want to make.

\begin{remark}\label{T} In \cite{DP}, section 6, 7 and 8 the authors use the fact
that when $|\nabla u|\geq c$, $u$ satisfies  a linear
nondivergence uniformly elliptic equation, $Tu=0$. In our case we
have that  when $|\nabla u|\geq c$, $u$ is a solution of the
equation defined in Remark \ref{lineal}. As in those sections the
authors only use the fact that this operator is linear and
uniformly elliptic, then the results of those sections in \cite{DP} extend to our
case without any change.
\end{remark}

For the reader convenience,  we will sketch here the proof of the regularity of the free
boundary by a series of steps  and we will
write down the proof in those cases in which we had to make  modifications.

\subsection{Flatness and nondegeneracy of the gradient}

\medskip

\begin{defi}[Flat free boundary points] Let $0<\sigma_+,
\sigma_-\le 1$ and $\tau>0$. We  say that $u$ is of class
$$
F(\sigma_+,\sigma_-; \tau)\quad \mbox{in}\quad B_\rho=B_\rho(0)
$$
if
\begin{enumerate}
\item $0\in \partial\{u>0\}$ and
$$
\begin{array}{ll}
u=0 & \mbox{for}\quad x_N\ge \sigma_+ \rho,\\
u(x)\ge -\lambda^*(x_N+\sigma_-\rho) & \mbox{for}\quad x_N\le
-\sigma_-\rho.
\end{array}
$$
\item $|\nabla u|\le \lambda^* (1+\tau)$ in $B_\rho$.
\end{enumerate}
If the origin is replaced by $x_0$ and the direction $e_N$ by the
unit vector $\nu$ we say that $u$ is of class
$F(\sigma_+,\sigma_-; \tau)$ in $B_\rho(x_0)$ in direction $\nu$.
\end{defi}

It is in the proof of the following theorems where we strongly use
the extra hypothesis (5) of weak solution (II) (which is
always satisfied by weak solutions (I)). For the details see
Section 6 in \cite{DP}.

\begin{teo}\label{flat2}
There exists $\sigma_0>0$ and $C_0>0$ such that
$$u\in F(\sigma,1; \sigma)\ in\ B_1 \ \ implies \ \ u \in
F(2\sigma,C_0\sigma;\sigma) \ in\ B_{1/2}$$ for
$0<\sigma<\sigma_0$.
\end{teo}
\begin{proof}
It follows as in the proof of Theorem 6.3 in \cite{DP} by Remark
\ref{T}.
\end{proof}
\begin{teo}\label{nodegegrad}
For every $\delta>0$ there exists $\sigma_{\delta}>0$ and
$C_{\delta}>0$ such that
$$u\in F(\sigma,1; \sigma)\ in\ B_1 \ \ implies \ \ |\nabla u|\geq \lambda^*-\delta
\ in\ B_{1/2}\cap \{x_N\leq -C_{\delta} \sigma\}$$ for
$0<\sigma<\sigma_{\delta}$.
\begin{proof}
It follows as in the proof of Theorem 6.4 in \cite{DP} by Remark
\ref{T}.
\end{proof}

\end{teo}

\medskip

\subsection{Nonhomogeneous blow-up}

\begin{lema}\label{9.1}
Let $u_k\in F(\sigma_k,\sigma_k;\tau_k) \in B_{\rho_k}$ with
$\sigma_k \to 0$, $\tau_k \sigma_k^{-2}\to 0$. For $y\in B_1'$,
set
\begin{align*}&f_k^+(y)= \sup\{h: (\rho_k y, \sigma_k \rho_kh)\in
\partial\{u_k>0\}\},\\& f_k^-(y)= \inf\{h: (\rho_k y, \sigma_k
\rho_kh)\in \partial\{u_k>0\}\}.\end{align*} Then, for a
subsequence,
\begin{enumerate}
\item
$f(y)=\limsup_{\stackrel{z\to y} {k\to \infty}}f_k^+(z)=
\liminf_{\stackrel{z\to y} {k\to \infty}}f_k^-(z)\ \ for\ all\
y\in B_1'.$

\smallskip

 Further, $f_k^+\to f$, $f_k^-\to f$ uniformly,
 $f(0)=0$, $|f|\leq 1$ and $f$ is continuous.

\medskip

 \item $f$ is
subharmonic.
\end{enumerate}
\end{lema}
\begin{proof}
(1) is the analogue of Lemma 5.3 in \cite{ACF}. The proof is based
on Theorem 6.3 and is identical to the one of Lemma 7.3 in
\cite{AC}.

The proof of (2) is a little bit different since here we don't
have in general that $q_{u_k}(x)=g(\lambda^*)$ $\H- a.e $ point in
$\partial\{u_k>0\}$. Instead, we have that this equality holds
for  $\H- a.e$ point in $\partial_{red}\{u_k>0\}$.

We may assume by replacing $u_k$ by
$\widetilde{u}_k=\frac{1}{\rho_k}u_k(\rho_k x)$, that $\rho_k=1$.
Let us assume, by contradiction that there is a ball
$B'_{\rho}(y_0)\subset B_1'$ and a harmonic function $g$ in a
neighborhood of this ball, such that
$$
g>f  \mbox{ on } \partial B'_{\rho}(y_0)\quad  \mbox{ and
} \quad f(y_0)>g(y_0).
$$
Let,
$$Z^{+} =\{x\in B_1\,/\, x=(y,h),\ y\in B'_{\rho}(y_0),
h>\sigma_k g(y)\},$$ and similarly $Z_0$ and $Z^-$. As in Lemma
7.5 in \cite{AC}, using  the same test function and the
Representation Theorem \ref{repteo} (see Remark \ref{repweak})  we
arrive at,
\begin{equation}\label{igualrep}\int_{\{u_k>0\}\cap Z_0} g(|\nabla u_k|) \frac{\nabla
u_k}{|\nabla u_k|} \cdot \nu\, d\H= \int_{\partial\{u_k>0\}\cap
Z^+} q_{u_k}(x) \, d\H.\end{equation}
As $u_k\in
F(\sigma_k,\sigma_k,\tau_k)$ we have that $|\nabla u_k|\leq
\lambda^*(1+\tau_k)$ and,  by Lemma \ref{glambda*}, there holds that
$q_{u_k}(x)=g(\lambda^*)$ for $\H- a.e$ point in
$\partial_{red}\{u_k>0\}$. Then, by \eqref{igualrep} we have,
\begin{equation}\label{desilambda1}g(\lambda^*) \H(\partial_{red}\{u_k>0\}\cap Z^+)\leq
g(\lambda^*(1+\tau_k)) \H(\{u_k>0\}\cap Z_0).\end{equation} On the
other hand, by the excess area estimate in Lemma 7.5 in \cite{AC}
we have that, $$\H(\partial_{red}E_k\cap Z)\geq
\H(Z_0)+c\sigma_k^2,$$ where $Z=B'_{\rho}(y_0)\times \mathbb{R}$
and $E_k=\{u_k>0\}\cup Z^-.$

We also have,
$$\H(\partial_{red}E_k\cap
Z)\leq \H(Z^+\cap\partial_{red}\{u_k>0\})+\H(Z_0\cap\{u_k=0\}).$$
Using these two inequalities and the fact that
$\H(Z_0\cap\partial\{u_k>0\})=0$ (if this is not true we replace $g$ by $g+c_0$ for a small
constant $c_0$) we have that,
\begin{equation}\label{exces}\H(\partial_{red}\{u_k>0\}\cap Z^+)\geq \H(Z_0\cap
\{u_k>0\} )+c\sigma_k^2.\end{equation} Finally by
\eqref{desilambda1} and \eqref{exces} we have that,
$$g(\lambda^*) \big[\H(\{u_k>0\}\cap Z_0)+c\sigma_k^2\big]\leq
g(\lambda^*(1+\tau_k)) \H(\{u_k>0\}\cap Z_0).$$ Therefore, for
some positive constant $c$ we have
$$c\leq \frac{g(\lambda^*(1+\tau_k))-g(\lambda^*)}{\sigma_k^2}$$
and this contradicts the fact that $\frac{\tau_k}{\sigma_k^2}\to
0$ as $k\to \infty$.
\end{proof}

\begin{lema}
There exists a positive constant $C=C(N)$ such that, for any $y\in
B'_{r/2}$,
$$\int_0^{1/4} \frac{1}{r^2} \Big(\pint_{\partial
B'_r(y)}f-f(y)\Big)\leq C_1.$$

\end{lema}
\begin{proof}
It follows as in Lemma 8.3 at \cite{DP} ,  by Remark \ref{T} and
Theorem \ref{nodegegrad}.
\end{proof}

With these two lemmas we have by Lemma 7.7 and Lemma 7.8 in
\cite{AC},

\begin{lema}\label{flat}\begin{enumerate}

\item $f$ is Lipschitz in $\bar{B}'_{1/4}$ with Lipschitz constant
depending on $C_1$ and $N$. \item There exists a constant
$C=C(N)>0$ and for $0<\theta<1$, there exists
$c_{\theta}=c(\theta,N)>0$, such that we can find a ball $B_r'$
and a vector $l\in \mathbb{R}^{N-1}$
 with
$$c_{\theta} \leq r\leq \theta,\ \ |l|\leq C, \ \ \mbox{ and }
f(y)\leq l .y+\frac{\theta}{2} r \ \ \mbox{ for } |y|\leq r.$$

 \end{enumerate}
\end{lema}
And as in Lemma 7.9 in \cite{AC} we have,
\begin{lema}\label{flat31}
Let $\theta$, $C$, $c_{\theta}$ as in Lemma \ref{flat}.  There
exists a positive constants $\sigma_{\theta}$,  such that
\begin{equation}\label{flat51}u\in F(\sigma,\sigma;\tau) \mbox{ in } B_{\rho} \mbox{ in
direction } \nu\end{equation} with $\sigma\leq \sigma_{\theta}, \
\tau\leq \sigma_{\theta} \sigma^2$, implies
$$u\in F(\theta\sigma,1; \tau) \mbox{ in } B_{\bar{\rho}} \mbox{ in direction }
\bar{\nu}$$ for some $\bar{\rho}$ and $\bar{\nu}$ with
$c_{\theta}\rho\leq \bar{\rho}\leq \theta \rho$ and
$|\bar{\nu}-\nu|\leq C\sigma$, where
$\sigma_{\theta}=\sigma_{\theta}(\theta,N)$.
\end{lema}

\begin{lema}\label{flat3}
Given $0<\theta<1$, there exist positive constants
$\sigma_{\theta}$,  $c_{\theta}$ and $C$ such that
\begin{equation}\label{flat5}u\in F(\sigma,1;\tau) \mbox{ in } B_{\rho} \mbox{ in
direction } \nu\end{equation} with $\sigma\leq \sigma_{\theta}
\mbox{ and } \ \tau\leq \sigma_{\theta} \sigma^2$, then
$$u\in F(\theta\sigma,\theta\sigma;\theta^2 \tau) \mbox{ in } B_{\bar{\rho}} \mbox{ in direction }
\bar{\nu}$$ for some $\bar{\rho}$ and $\bar{\nu}$ with
$c_{\theta}\rho\leq \bar{\rho}\leq \frac{1}{4} \rho$ and
$|\bar{\nu}-\nu|\leq C\sigma$, where
$c_{\theta}=c_{\theta}(\theta,N)$, $C=C(N,\delta,g_0)$,
$\sigma_{\theta}=\sigma_{\theta}(\theta,N)$.
\end{lema}
\begin{proof}
We obtain the improvement of the value $\tau$ inductively. Assume
that $\rho=1$. If $\sigma_{\theta}$ is small enough, we can apply
Theorem \ref{flat2} and obtain
$$u\in F(C\sigma,C\sigma;\tau) \mbox{ in } B_{1/2} \mbox{ in direction } \nu.$$ Then
for $0<\theta_1\leq \frac{1}{2}$ we can apply Lemma \ref{flat31},
if again $\sigma_{\theta}$ is small, and we obtain
\begin{equation}\label{theta1}u\in F(C\theta_1\sigma,C\sigma;\tau) \mbox{ in }
B_{r_1} \mbox{ in direction } \nu_1\end{equation} for some $r_1,
\nu_1$ with
$$c_{\theta_1}\leq 2 r_1\leq \theta_1, \mbox{ and }
|\nu_1-\nu|\leq C\sigma.
$$

In order to improve $\tau$, we consider the
functions $U_{\ep}=\big(G(|\nabla u|) - G(\lambda^*) - \ep\big)^+$
and $U_0=\big(G(|\nabla u|) - G(\lambda^*)\big)^+$ in $B_{2r_1}$.
By Lemma \ref{const}, and (4) in Definitions \ref{weak} and
\ref{weak2} we know that $U_{\ep}$ vanishes in a neighborhood of
the free boundary. Since $U_{\ep}>0$ implies $G(|\nabla
u|)>G(\lambda^*)+\ep$, the closure of $\{U_{\ep}>0\}$ is contained
in $\{G(|\nabla u|)>G(\lambda^*)+\ep/2\}$.
 The function $u$ satisfies the linearized equation
 $$Tu=b_{ij}(\nabla u) D_{i j}u=0$$ where
 $b_{ij}$ is defined in \eqref{opT}, and is
  uniformly elliptic in $\{G(|\nabla
u|)>G(\lambda^*)+\ep/2\}$ with ellipticity constant $\beta$ independent of $u$.

Let $v=G(|\nabla u|)$. By Lemma 1 in \cite{Li2}, we have that $v$
satisfies,
$$Mv=D_i(b_{ij}(\nabla u) D_j v)\geq 0 \ \ \mbox{ in
}\{G(|\nabla u|)>G(\lambda^*)+\ep/2\}.$$ Hence $U_{\ep}$ satisfies
$$MU_{\ep}\geq 0 \ \ \mbox{ in } \{G(|\nabla u|)>G(\lambda^*)+\ep/2\}.$$
Extending the operator $M$ with the uniformly elliptic
divergence-form operator
$$\widetilde{M}w=D_i(\widetilde{b}_{ij}(x) D_j w)\ \ \mbox{ in }
B_{2r_1}$$ with measurable coefficients such that
$$\widetilde{b}_{ij}(x)= b_{ij}(\nabla u)
 \ \ \mbox{ in } \{G(|\nabla u|)>G(\lambda^*)+\ep/2\},$$
we obtain
$$\widetilde{M}U_{\ep}\geq 0\ \ \mbox{ in } B_{2r_1}.$$
Moreover, by \eqref{flat5} we have that  $U_{\ep}\leq
G(\lambda^*(1+\tau))-G(\lambda^*)$ and
 by \eqref{theta1}
$U_{\ep}=0$ in $B=B_{r_1/4}\big(\frac{r_1}{2}\nu_1\big)$, if
$C\sigma\leq 1/2$.

Take now, $V$ such that, $$\begin{cases} \widetilde{M}V=0 & \mbox{
in } B_{2r_1}\setminus \bar{B},\\
V=G(\lambda^*(1+\tau))-G(\lambda^*) & \mbox{ on } \partial
B_{2r_1},\\
V=0 &\mbox{ on } \partial B.
\end{cases}$$
Then, there exists $0<c(N,\beta)<1$ such that $V\leq c
\big(G(\lambda^*(1+\tau))-G(\lambda^*)\big)$ in $B_{r_1}$. Applying the
maximum principle we have that, $U_{\ep}\leq c
(G(\lambda^*(1+\tau))-G(\lambda^*))$ in $B_{r_1}$. Taking $\ep\to
0$ we obtain,
$$\begin{aligned}G(|\nabla u|)&\leq c
G(\lambda^*(1+\tau))+G(\lambda^*)(1-c)\quad \mbox{ in } B_{r_1}.
\end{aligned}
$$
Since, $G(\lambda^*(1+\tau))=G(\lambda^*)+g(\lambda^*)
\lambda^*\tau+o(\tau)$ we have that
$$c
G(\lambda^*(1+\tau))+G(\lambda^*)(1-c)=G(\lambda^*)+c g(\lambda^*)
\lambda^*\tau+o(\tau),$$ and since $G$ is strictly increasing, we
have, $$\begin{aligned}|\nabla u| &\leq G^{-1}(G(\lambda^*)+c
g(\lambda^*)
\lambda^*\tau+o(\tau))\\&=\lambda^*+\frac{1}{g(\lambda^*)}
(g(\lambda^*) \lambda^* \tau c +o(\tau))+ o(\tau)\\&=
\lambda^*\Big(1+
 \tau \big(c +\frac{o(\tau)}{\tau}\big)\Big)\leq
\lambda^* \Big(1+
 \tau \frac{(c+1)}{2}\Big),
\end{aligned}$$
if we choose $\tau$ small enough.
 And we see that if we choose $\theta_1$ small enough (depending
on $N$), we have
$$u\in F(\theta_0\sigma,1;\theta_0^2\tau) \mbox{ in } B_{r_1} \mbox{ in
direction } \nu_1,$$ where $\theta_0=\sqrt{\frac{c+1}{2}}$.

We can repeat this argument a finite number of times, and we obtain
$$u\in F(\theta_0^m\sigma,1;\theta_0^{2m}\tau) \mbox{ in } B_{r_1...r_m} \mbox{ in
direction } \nu_m,$$ with
$$c_{\theta_j}\leq 2 r_j\leq \theta_j, \mbox{ and }
|\nu_m-\nu|\leq \frac{C}{1-\theta_0}\sigma.$$ Finally we choose
$m$ large enough and use Theorem \ref{flat2}.
\end{proof}

\subsection{Smoothness of the free boundary}

\begin{teo}\label{regfin}
Suppose that $u$ is a weak solution, and $D\subset\subset\Omega$.
Then there exist positive constants $\bar{\sigma}_0$, $C$ and
$\alpha$ such that if
$$u\in F(\sigma,1;\infty) \quad in\ B_{\rho}(x_0)\subset D \mbox{
in direction } \nu$$ with $\sigma\leq \bar{\sigma}_0$, $\rho\leq
\bar{\rho}_0(\bar{\sigma}_0,\sigma)$, then
$$B_{\rho/4}(x_0)\cap \partial\{u>0\} \mbox{ is a } C^{1,\alpha}
\mbox{ surface,}$$ more precisely, a graph in direction $\nu$ of a
$C^{1,\alpha}$ function, and, for any $x_1$, $x_2$ on this surface
$$|\nu(x_1)-\nu(x_2)|\leq C\sigma
\Big|\frac{x_1-x_2}{\rho}\Big|^{\alpha}$$
\end{teo}
\begin{proof}
By property (4) in Definitions \ref{weak} and \ref{weak2} we have that,
 for every   $\rho-$ neighborhood  $D_{\rho}$ of
$D\cap\partial\{u>0\}$,
\begin{align*}
 |\nabla u(x)| \leq
\lambda^*+\tau(\rho),\qquad \mbox{for every } x\in D_{\rho}
\end{align*}
where $\tau(\rho)\rightarrow 0$ when $\rho\to 0$.

Therefore, $$u\in F(\sigma,1;\tau) \quad in\ B_{\rho}(x_0) \mbox{
in direction } \nu.$$

Applying Theorem \ref{flat2} we have that $$u\in
F(C_0\sigma,C_0\sigma;\tau) \quad in\ B_{\rho/2}(x_0) \mbox{ in
direction } \nu$$ if $\sigma\leq\sigma_0$ and $\tau\leq \sigma$.

Let $x_1\in B_{\rho/2}(x_0) \cap\partial\{u>0\}$ then
$$u\in
F(C_0\sigma,1;\tau) \quad in\ B_{\rho/2}(x_1) \mbox{ in direction
} \nu$$ and applying again Theorem \ref{flat2} we have,
$$u\in
F(C_0^2\sigma,C_0^2\sigma;\tau) \quad in\ B_{\rho/4}(x_1) \mbox{
in direction } \nu$$ if $C_0 \sigma \leq \sigma_0$ and $\tau\leq
C_0 \sigma$.

Let $0<\theta<1$, take $\rho_0=\rho/4$, $\nu_0=\nu$, $C=C_0^2$,
$\sigma\leq\frac{\sigma_{\theta}}{C}$ and $\tau\leq
\sigma_{\theta} C^2\sigma^2$. Now, by Lemma \ref{flat3} and
iterating we get that there exist  sequences  $\rho_m$ and $\nu_m$
such that,
$$u\in
F(\theta^mC\sigma,\theta^m C\sigma;\theta^{2m}\tau) \quad in\
B_{\rho_m}(x_1) \mbox{ in direction } \nu_m$$ with $c_{\theta}
\rho_m\leq \rho_{m+1}\leq \rho_m/4$ and $|\nu_{m+1}-\nu_m|\leq
\theta^m C \sigma$.

Thus, we have that $|\langle x-x_1,\nu_m\rangle|\leq \theta^m C\sigma
\rho_m$ for $x\in B_{\rho_m}(x_1)\cap\partial\{u>0\}$.

We also have that there exists $\nu(x_1)=\lim_{m\to\infty} \nu_m$  and
$$|\nu(x_1)-\nu_m|\leq \frac{C\theta^m}{1-\theta} \sigma.$$

Now let $x\in B_{\rho/4}(x_1)\cap\partial\{u>0\}$ and choose $m$
such that $\rho_{m+1}\leq |x-x_1|\leq \rho_m$. Then
$$|\langle x-x_1, \nu(x_1)\rangle|\leq C \theta^m
\sigma\Big(\frac{|x-x_1|}{1-\theta}+\rho_m\Big) \leq C \theta^m
\sigma \Big(\frac{1}{1-\theta}+\frac{1}{c_{\theta}}\Big)|x-x_1|
$$
and since $|x-x_1|\leq c_{\theta}^{m+1} \rho_0$ we have
$$\theta^{m+1}\leq \Big(\frac{|x-x_1|}{\rho_0}\Big)^{\alpha}\quad \mbox{with }
\alpha=\frac{log(\theta)}{log(c_{\theta})},$$ and we conclude that
$$|\langle x-x_1, \nu(x_1)\rangle|\leq
\frac{C\sigma}{\rho^{\alpha}}|x-x_1|^{1+\alpha}.$$

 Finally, observe that the result follows if we take,
$\bar{\sigma}_0=\min\{\sigma_0,\frac{\sigma_0}{C_0},\frac{\sigma_{\theta}}{C}\}$
and if we choose $\bar{\rho}_0$ small enough such that if
$\rho\leq \bar{\rho}_0$,  $\tau(\rho)\leq \min\{\sigma, C_0
\sigma, {\sigma}_{\theta} C^2 \sigma^2\}$.
\end{proof}

\begin{remark}\label{remfin} By Lemma \ref{blow4}, Definition \ref{weak2} and by the nondegeneracy,
 we have that
there exists a set $A\subset \partial_{red}\{u>0\}$, with
$\H(\partial_{red}\{u>0\}\setminus A)=0$, such that for $x_0\in A$
we have that
 $u\in F(\sigma_{\rho},1;\infty)$ in $B_{\rho}(x_0)$ in direction
$\nu_u(x_0)$, with $\sigma_{\rho}\to 0$ for $\rho\to 0$. Observe
that by Theorem \ref{blow2} when $u$ is a minimizer
$A=\partial_{red}\{u>0\}$. Hence applying Theorem \ref{regfin} we
have,\end{remark}

\begin{teo}\label{teo.regularity}
If $u$ is a  weak solution then there exists a subset $A\subset\partial_{red}\{u>0\}$ with
$\H (\partial_{red}\{u>0\}\setminus A)=0$ such that for any $x_0\in A$ there exists $r>0$ so that
$B_r(x_0)\cap\partial\{u>0\}$ is a $C^{1,\alpha}$ surface.
Moreover, if $u$ satisfies Definition \ref{weak} then the
remainder of $\partial\{u>0\}$ has $\H$--measure zero. Finally, if $u$ is a minimizer, $\partial_{red}\{u>0\}$
is a $C^{1,\alpha}$ surface and $\H(\partial\{u>0\}\setminus \partial_{red}\{u>0\}=0$.
\end{teo}

\end{document}